\magnification=\magstep1
\hsize=6.5truein
\hoffset=0.0truein
\baselineskip 1.3\normalbaselineskip

\tolerance=10000
\def\sqr{$\vcenter{\hrule height .3mm
\hbox {\vrule width .3mm height 2mm \kern 2mm
\vrule width .3mm} \hrule height .3mm}$}

\input amssym.def
\input amssym.tex

\def \Proof{\noindent {\bf Proof.\ \ }}

\def \R{{\Bbb R}}
\def \E{{\Bbb E}}
\def \P{{\Bbb P}}

\def \N{{\Bbb N}}

\def \seq#1#2{#1_1,\dots,#1_#2}

\def \sm#1#2{\sum_{#1=1}^#2}
\def \ge{\geqslant}
\def \le{\leqslant}

\def \e{\epsilon}

\def \g{\gamma}
\def \d{\delta}

\def \a{\alpha}
\def \b{\beta}

\def \sp#1{\langle#1\rangle}
\def \ra{\rightarrow}
 
\def \cb {{\cal B}}
\def \ch {{\cal H}}
\def \cj {{\cal J}}
\def \cd {{\cal D}}
\def \ck {{\cal K}}

\def \cp {{\cal P}}
\def \cq {{\cal Q}}
\def \Hom {\mathop{\rm Hom}}

\def \Z{{\Bbb Z}}
\def \oct{\mathop{{\rm Oct}}}

\def \var{\mathop{\rm var}}

\centerline {\bf Hypergraph Regularity and the multidimensional
Szemer\'edi Theorem.}
\medskip

\centerline {\bf W. T. Gowers}
\bigskip

\noindent {\bf Abstract.} {\sl We prove analogues for hypergraphs of 
Szemer\'edi's regularity lemma and the associated counting lemma
for graphs. As an application, we give the first combinatorial
proof of the multidimensional Szemer\'edi theorem of Furstenberg
and Katznelson, and the first proof that provides an explicit bound.
Similar results with the same consequences have been obtained 
independently by Nagle, R\"odl, Schacht and Skokan.}
\bigskip

\noindent {\bf \S 1. Introduction.}
\medskip

Szemer\'edi's theorem states that, for every real number $\d>0$ and
every positive integer $k$, there exists a positive integer $N$ such
that every subset $A$ of the set $\{1,2,\dots,N\}$ of size at least
$\d N$ contains an arithmetic progression of length $k$. There are now
three substantially different proofs of the theorem, Szemer\'edi's
original combinatorial argument [Sz1], an ergodic-theory proof due to
Furstenberg (see for example [FKO]) and a proof by the author using
Fourier analysis [G1].  Interestingly, there has for some years been a
highly promising programme for yet another proof of the theorem,
pioneered by Vojta R\"odl (see for example [R]), developing an argument 
of Ruzsa and Szemer\'edi [RS] that proves the result for progressions 
of length three. Let us briefly sketch their argument.

The first step is the famous regularity lemma of Szemer\'edi [Sz2]. If 
$G$ is a graph and $A$ and $B$ are sets of vertices in $V$, then
let $e(A,B)$ stand for the number of pairs $(x,y)\in A\times B$
such that $xy$ is an edge of $G$. Then the {\it density} $d(A,B)$ 
of the pair $(A,B)$ is $e(A,B)/|A||B|$. The pair is $\e$-{\it regular} 
if $|d(A',B')-d(A,B)|\le\e$ for all subsets $A'\subset A$ and 
$B'\subset B$ such that $|A'|\ge\e|A|$ and $|B'|\ge\e|B|$. The basic
idea is that a pair is regular with density $d$ if it resembles a 
random graph with edge-probability $d$. Very roughly, the regularity 
lemma asserts that every graph can be decomposed into a few pieces,
almost all of which are random-like. The precise statement is as
follows.

\proclaim Theorem {1.1}. Let $\e>0$. Then there exists a positive
integer $K_0$ such that, given any graph $G$, the vertices can be
partitioned into $K\leq K_0$ sets $V_i$, with sizes differing by at
most 1, such that all but at most $\e K^2$ of the pairs $(V_i,V_j)$
are $\e$-regular.

\noindent A partition is called $\e$-{\it regular} if it satisfies
the conclusion of Theorem 1.1. (Note that we allow $i$ to equal $j$ in
the definition of a regular pair, though if $K$ is large then this
does not make too much difference.) The regularity lemma is
particularly useful in conjunction with a further result, known as the
counting lemma. To state it, it is very convenient to use the notion
of a graph homomorphism. If $G$ and $H$ are graphs, then a function
$\phi:V(H)\ra V(G)$ is called a {\it homomorphism} if $\phi(x)\phi(y)$
is an edge of $G$ whenever $xy$ is an edge of $H$. It is an {\it
isomorphic embedding} if in addition $\phi(x)\phi(y)$ is not an edge
of $G$ whenever $xy$ is not an edge of $H$.

\proclaim Theorem {1.2}. For every $\a>0$ and every $k$ there exists
$\e>0$ with the following property. Let $\seq V k$ be sets of vertices 
in a graph $G$, and suppose that for each pair $(i,j)$ the pair $(V_i,V_j)$
is $\e$-regular with density $d_{ij}$. Let $H$ be a graph with vertex
set $(\seq x k)$, let $v_i\in V_i$ be chosen independently and uniformly 
at random, and let $\phi$ be the map that takes $x_i$ to $v_i$ for each
$i$. Then the probability that $\phi$ is an isomorphic embedding differs
from $\prod_{x_ix_j\in H}d_{ij}\prod_{x_ix_j\notin H}(1-d_{ij})$ by
at most $\a$. 

\noindent Roughly, this result tells us that the $k$-partite graph
induced by the sets $\seq V k$ contains the right number of labelled
induced copies of the graph $H$. Let us briefly see why this result is true
when $H$ is a triangle. Suppose that $U,V,W$ are three sets of
vertices and the pairs $(U,V)$, $(V,W)$ and $(W,U)$ are $\e$-regular
with densities $\zeta$, $\eta$ and $\theta$ respectively. Then a
typical vertex of $U$ has about $\zeta|V|$ neighbours in $V$ and
$\theta|W|$ neighbours in $W$. By the regularity of the pair $(V,W)$,
these two neighbourhoods span about $\eta(\zeta|V|)(\theta|W|)$ edges
in $G$, creating that many triangles. Summing over all vertices of $U$ 
we obtain the result.

The next step in the chain of reasoning is the following 
innocent-looking statement about graphs with few triangles. Some of the
details of the proof will be sketched rather than given in full.

\proclaim Theorem {1.3}. For every constant $a>0$ there exists a
constant $c>0$ with the following property. If $G$ is any graph with 
$n$ vertices that contains at most $cn^3$ triangles, then it is
possible to remove at most $an^2$ edges from $G$ to make it triangle-free.

\Proof This theorem is a simple consequence of the regularity lemma. 
Indeed, let $\e=\e(a)>0$ be sufficiently small and let $\seq V K$ be 
an $\e$-regular partition of the vertices of $G$.  If there are fewer
than $a|V_i||V_j|/100$ edges between $V_i$ and $V_j$, then remove all
those edges, and also remove all edges from $V_i$ to $V_j$ if
$(V_i,V_j)$ is not an $\e$-regular pair. Since the partition is
$\e$-regular, we have removed fewer than $an^2$ edges, and the
resulting graph must either be triangle-free or contain several
triangles. To see why this is, suppose that $(x,y,z)$ is a triangle in
$G$ (after the edges have been removed), and suppose that $(x,y,z)\in
V_i\times V_j\times V_k$. Then by our construction the pair
$(V_i,V_j)$ must be regular and must span many edges (because we did
not remove the edge $(x,y)$) and similarly for the pairs $(V_j,V_k)$
and $(V_i,V_k)$. But then, by the counting lemma for triangles, the
sets $V_i$, $V_j$ and $V_k$ span at least $a^3|V_i||V_j||V_k|/10^6$
triangles. Each $V_i$ has cardinality at least $n/2K$, where $K$ 
depends on $\e$ only (which itself depends on $a$ only). This proves 
that the result is true provided that $c\le a^3/2^310^6K^3$. \hfill
$\square$ \bigskip

Ruzsa and Szemer\'edi [RS] observed that Theorem 1.3 implies Szemer\'edi's
theorem for progressions of length 3. More recently, Solymosi noticed
[So1,2] that it also implied the following two-dimensional
generalization. (Actually, neither of these statements is quite
accurate. There are several closely related graph-theoretic results
that have these consequences and can be proved using the regularity
lemma, of which Theorem 1.3 is one. Ruzsa and Szemer\'edi and
Solymosi did not use Theorem 1.3 itself but their arguments are not
importantly different.)

\proclaim Corollary {1.4}. For every $\d>0$ there exists $N$ such
that every subset $A\subset[N]^2$ of size at least $\d N^2$ contains
a triple of the form $(x,y)$, $(x+d,y)$, $(x,y+d)$ with $d>0$. 

\Proof  First, note that an easy argument allows us to replace
$A$ by a set $B$ that is symmetric about some point. Briefly, if
the point $(x,y)$ is chosen at random then the intersection of $A$ 
with $(x,y)-A$ has expected size $c\d^2N^2$ for some absolute 
constant $c>0$, lives inside the grid $[-N,N]^2$, and has the
property that $B=(x,y)-B$. So $B$ is still reasonably dense, and
if it contains a subset $K$ then it also contains a translate of
$-K$. So we shall not worry about the condition $d>0$. (I am grateful
to Ben Green for bringing this trick to my attention. As it happens,
the resulting improvement to the theorem is something of a side issue, 
since the positivity of $d$ does not tend to be used in applications. 
See for instance Corollary 1.5 below. See also the remark at the
beginning of the proof of Theorem 10.3.) 

Without loss of generality, the original set $A$ is symmetric in
this sense. Let $X$ be the set of all vertical lines through $[N]^2$,
that is, subsets of the form $\{(x,y):x=u\}$ for some $u\in[N]$.
Similarly, let $Y$ be the set of all horizontal lines. Define
a third set, $Z$, of diagonal lines, that is, lines of constant
$x+y$. These sets form the vertex sets of a tripartite graph, 
where a line in one set is joined to a line in another if and only
if their intersection belongs to $A$. For example, the line
$x=u$ is joined to the line $y=v$ if and only if $(u,v)\in A$
and the line $x=u$ is joined to the line $x+y=w$ if and only 
if $(u,w-u)\in A$. 

Suppose that the resulting graph $G$ contains a triangle of lines
$x=u$, $y=v$, $x+y=w$. Then the points $(u,v)$, $(u,w-u)$
and $(w-v,v)$ all lie in $A$. Setting $d=w-u-v$, we can rewrite
them as $(u,v)$, $(u,v+d)$, $(u+d,v)$, which shows that we are
done unless $d=0$. When $d=0$, we have $u+v=w$, which corresponds
to the degenerate case when the vertices of the triangle in $G$
are three lines that intersect in a single point. Clearly, this
can happen in at most $|A|=o(N^3)$ ways.

Therefore, if $A$ contains no configuration of the desired kind,
then the hypothesis of Theorem 1.3 holds, and we can remove 
$o(N^2)$ edges from $G$ to make it triangle-free. But this is
a contradiction, because there are at least $\d N^2$ degenerate 
triangles and they are edge-disjoint. \hfill $\square$ \bigskip

An easy consequence of Corollary 1.4 is the case $k=3$ of Szemer\'edi's 
theorem, which was first proved by Roth [R] using Fourier analysis.

\proclaim Corollary {1.5}. For every $\d>0$ there exists $N$ such
that every subset $A$ of $\{1,2,\dots,N\}$ of size at least $\d N$
contains an arithmetic progression of length 3.

\Proof  Define $B\subset[N]^2$ to be the set of all $(x,y)$ such
that $x+2y\in A$. It is straightforward to show that $B$ has 
density at least $\eta>0$ for some $\eta$ that depends on $\d$
only. Applying Corollary 1.2 to $B$ we obtain inside it three points
$(x,y)$, $(x+d,y)$ and $(x,y+d)$. Then the three numbers
$x+2y$, $x+d+2y$ and $x+2(y+d)$ belong to $A$ and form an
arithmetic progression. \hfill $\square$ \bigskip

And now the programme for proving Szemer\'edi's theorem in general
starts to become clear. Suppose, for example, that one would like to
prove it for progressions of length 4. After a little thought, one
sees that the direction in which one should generalize Theorem 1.3 is
the one that takes graphs to 3-{\it uniform hypergraphs}, or 3-{\it
graphs}, for short, which are set systems consisting of subsets of
size 3 of a set $X$ (just as a graph consists of pairs). If $H$ is a
3-uniform hypergraph, then a {\it simplex} in $H$ is a set of four
vertices $x,y,z$ and $w$ of $H$ (that is, elements of the set $X$)
such that the four triples $xyz$, $xyw$, $xzw$ and $yzw$ all belong to
$H$. The following theorem of Frankl and R\"odl is a direct
generalization of Theorem 1.3, but its proof is much harder.

\proclaim Theorem {1.6}. For every constant $a>0$ there exists a 
constant $c>0$ with the following property. If $H$ is any 3-uniform 
hypergraph with $n$ vertices that contains at most $cn^4$ simplices,
then it is possible to remove at most $an^3$ edges from $H$ to make
it simplex-free.

\noindent As observed by Solymosi, it is straightforward to generalize 
the proof of Theorem {1.4} and show that Theorem 1.6 has the following
consequence.

\proclaim Theorem {1.7}. For every $\d>0$ there exists $N$ such that
every subset $A\subset[N]^3$ of size at least $\d N^3$ contains a
quadruple of points of the form 
$$\{(x,y,z),(x+d,y,z),(x,y+d,z),(x,y,z+d)\}$$ 
with $d>0$. 

\noindent Similarly, Szemer\'edi's theorem for progressions of length 
four is an easy consequence of Theorem 1.7 (and once again one does
not need the positivity of $d$).

It may look as though this section contains enough hints to enable any
sufficiently diligent mathematician to complete a proof of the entire
theorem. Indeed, here is a sketch for the 3-uniform case. First, one
proves the appropriate 3-graph analogue of Szemer\'edi's regularity
lemma. Then, given a hypergraph $H$, one applies this lemma. Next, one
removes all sparse triples and all triples that fail to be regular. If
the resulting hypergraph contains a simplex, then any three of the
four sets in which its vertices lie must form a dense regular triple,
and therefore (by regularity) the hypergraph contains many simplices,
contradicting the original assumption.

The trouble with the above paragraph is that it leaves unspecified
what it means for a triple to be regular. It turns out to be
surprisingly hard to come up with an appropriate definition, where
``appropriate'' means that it must satisfy two conditions. First, it
should be weak enough for a regularity lemma to hold: that is, one
should always be able to divide a hypergraph up into regular
pieces. Second, it should be strong enough to yield the conclusion
that four sets of vertices, any three of which form a dense regular
triple, should span many simplices. The definition that Frankl and
R\"odl used for this purpose is complicated and it proved very hard to
generalize. In [G2] we gave a different proof which is in some ways
more natural. The purpose of this paper is to generalize the results
of [G2] from 3-uniform hypergraphs to $k$-uniform hypergraphs for
arbitrary $k$, thereby proving the full multidimensional version of
Szemer\'edi's theorem (Theorem 10.3 below), which was first proved by
Furstenberg and Katznelson [FK]. This is the first proof of the
multidimensional Szemer\'edi theorem that is not based on
Furstenberg's ergodic-theoretic approach, and also the first proof
that gives an explicit bound. The bound, however, is very weak---it
gives an Ackermann-type dependence on the initial parameters.

Although this paper is self-contained, we recommend reading [G2]
first. The case $k=3$ contains nearly all the essential ideas, and
they are easier to understand when definitions and proofs can be given
directly. Here, because we are dealing with a general $k$, many of the
definitions have to be presented inductively. The resulting proofs can
be neater, but they may appear less motivated if one has not examined
smaller special cases. For this reason, we do indeed discuss a special
case in the next section, but not in as complete a way as can be
found in [G2]. Furthermore, the bulk of [G2] consists of background
material and general discussion (such as, for example, a complete
proof of the regularity lemma for graphs and a detailed explanation of
how the ideas relate to those of the analytic approach to
Szemer\'edi's theorem in [G1]). Rather than repeat all that motivating
material, we refer the reader to that paper for it.

The main results of this paper have been obtained independently by
Nagle, R\"odl, Schacht and Skokan [NRS,RS]. They too prove hypergraph
generalizations of the regularity and counting lemmas that imply
Theorem 10.3 and Szemer\'edi's theorem. However, they formulate their
generalizations differently and there are substantial differences
between their proof and ours. Broadly speaking, they take the proof of
Frankl and R\"odl as their starting point, whereas we start with the
arguments of [G2]. This point is discussed in more detail in the 
introduction to \S 6 of this paper, and also at the end of [G2].
\bigskip

\noindent {\bf \S 2. A discussion of a small example.}
\medskip

The hardest part of this paper will be the proof of a counting lemma,
which asserts that, under certain conditions, a certain type of
structure ``behaves randomly'' in the sense that it contains roughly
the expected number (asymptotically speaking) of configurations of any
fixed size. In order even to state the lemma, we shall have to develop
quite a lot of terminology, and the proof will involve a rather
convoluted inductive argument with a somewhat strange inductive
hypothesis. The purpose of this section is to give some of the
argument in a special case. The example we have chosen is small
enough that we can discuss it without the help of the terminology
we use later: we hope that as a result the terminology will be
much easier to remember and understand (since it can be related
to the concrete example). Similarly, it should be much clearer
why the inductive argument takes the form it does. From a logical
point of view this section is not needed: the reader who likes
to think formally and abstractly can skip it and move to the next 
section\footnote{$^1$}{This section was not part of the original 
submitted draft. One of the referees suggested treating a small
case first, and when I reread the paper after a longish interval
I could see just how much easier it would be to understand if I
followed the suggestion}.

To put all this slightly differently, the argument is of the following
kind: there are some simple techniques that can be used quite
straightforwardly to prove the counting lemma in any particular
case. However, as the case gets larger, the expressions that appear
become quite long (as will already be apparent in the example we are
about to discuss), even if the method for dealing with them is
straightforward. In order to discuss the general case, one is forced
to {\it describe} in general terms what it is one is doing, rather
than just going ahead and doing it, and for that it is essential to
devise a suitably compact notation, as well as an inductive hypothesis
that is sufficiently general to cover all intermediate stages in the
calculation.

Now we are ready to turn to the example itself.
Let $X$, $Y$, $Z$ and $T$ be four finite sets. We shall adopt
the convention that variables that use a lower-case letter of
the alphabet range over the set denoted by the corresponding
upper-case letter. So, for example, $x'$ would range over $X$.
Similarly, if we refer to ``the function $v(y,z,t)$,'' it should
be understood that $v$ is a function defined on $Y\times Z\times T$.

For this example, we shall look at three functions, $f(x,y,z)$,
$u(x,y,t)$ and $v(y,z,t)$. (The slightly odd choices
of letters are deliberate: $f$ plays a different role from the
other functions and $t$ plays a different role from the other
variables.) We shall also assume that they are supported in a
quadripartite graph $G$, with vertex sets $X$, $Y$, $Z$ and $T$,
in the sense that $f(x,y,z)$ is non-zero only if $xy$, $yz$ and
$xz$ are all edges of $G$, and similarly for the other three
functions. As usual, we shall feel free to identify $G$ with
its own characteristic function, so another way of stating our
assumption is that $f(x,y,z)=f(x,y,z)G(x,y)G(y,z)G(x,z)$. 

We will need one useful piece of shorthand as the proof proceeds.
We shall write $f_{x,x'}(y,z)$ for $f(x,y,z)f(x',y,z)$, and
similarly for the other functions (including $G$) and variables.
We shall even iterate this, so that $f_{x,x',y,y'}(z)$ means
$$f(x,y,z)f(x',y,z)f(x,y',z)f(x',y',z).$$
Of particular importance to us will be the quantity 
$\oct(f)=\E_{x,x',y,y',z,z'}f_{x,x',y,y',z,z'}$, which is a 
count of octahedra, each one weighted by the product of
the values that $f$ takes on its eight faces. 

Now let us try to obtain an upper bound for the quantity
$$\E_{x,y,z,t}f(x,y,z)u(x,y,t)v(y,z,t).$$
Our eventual aim will be to show that this is small if
$\oct(f)$ is small and the six parts of $G$ are sufficiently 
quasirandom. However, an important technical idea of the proof, 
which simplifies it considerably, is to avoid using the quasirandomness
of $G$ for as long as possible. Instead, we make no assumptions
about $G$ (though we imagine it as fairly sparse and very 
quasirandom), and try to obtain an upper bound for our expression
in terms of $f_{x,x',y,y',z,z'}$ and $G$. Only later do we use the fact
that we can handle quasirandom graphs. In the more general 
situation, something similar occurs: now $G$ becomes a hypergraph,
but in a certain sense it is less complex than the original 
hypergraph, which means that its good behaviour can be assumed  
as the complicated inductive hypothesis alluded to earlier. 

As with many proofs in arithmetic combinatorics, the upper bound we
are looking for is obtained by repeated use of the Cauchy-Schwarz
inequality, together with even more elementary tricks such as
interchanging the order of expectation, expanding out the square of an
expectation, or using the inequality
$\E_xf(x)g(x)\leq\|f\|_1\|g\|_\infty$. The one thing that makes the
argument slightly (but only slightly) harder than several other
arguments of this type is that it is essential to use the
Cauchy-Schwarz inequality efficiently, and easy not to do so if one is
careless. In many arguments it is enough to use the inequality $(\E_x
f(x))^2\leq\E_x f(x)^2$, but for us this will usually be inefficient
because it will usually be possible to identify a small set of $x$
outside which $f(x)$ is zero. Letting $A$ be the characteristic
function of that set, we can write $f=Af$, and we then have the
stronger inequality $(\E_x f(x))^2\leq\E_xA(x)\E_x f(x)^2$.

Here, then, is the first part of the calculation that gives us
the desired upper bound. We need one further assumption: that the 
functions $f$, $u$ and $v$ take values in the interval $[-1,1]$.
$$\eqalign{\Bigl(&\E_{x,y,z,t}f(x,y,z)u(x,y,t)v(y,z,t)\Bigr)^8\cr
&=\Bigl(\E_{y,z,t}\E_xf(x,y,z)u(x,y,t)v(y,z,t)\Bigr)^8\cr
&=\Bigl(\E_{y,z,t}G(y,z)G(y,t)G(z,t)\E_xf(x,y,z)u(x,y,t)v(y,z,t)\Bigr)^8\cr
&\le\Bigl(\E_{y,z,t}G(y,z)G(y,t)G(z,t)\Bigr)^4
\Bigl(\E_{y,z,t}\Bigl(\E_xf(x,y,z)u(x,y,t)v(y,z,t)
\Bigr)^2\Bigr)^4.\cr}$$ The inequality here is Cauchy-Schwarz, and we
have used the fact that $v(y,z,t)$ is non-zero only if
$G(y,z)G(y,t)G(z,t)=1$.  For the same reason, the second bracket is at
most
$$\eqalign{\Bigl(\E_{y,z,t}\Bigl(\E_xf(x,y,z)u(x,y,t)
&G(y,z)G(y,t)G(z,t)\Bigr)^2\Bigr)^4\cr
&=\Bigl(\E_{y,z,t}\Bigl(\E_xf(x,y,z)u(x,y,t)
G(z,t)\Bigr)^2\Bigr)^4\cr
&=\Bigl(\E_{x,x'}\E_{y,z,t}f_{x,x'}(y,z)u_{x,x'}(y,t)
G(z,t)\Bigr)^4\cr
&\le\E_{x,x'}\Bigl(\E_{y,z,t}f_{x,x'}(y,z)u_{x,x'}(y,t)
G(z,t)\Bigr)^4\cr}$$
The first equality here follows from the fact that $G(y,z)$
and $G(y,t)$ are 1 whenever $f(x,y,z)$ and $u(x,y,t)$ are
non-zero. The inequality is a simple case of Cauchy-Schwarz,
applied twice.

Simple manipulations and arguments of the above kind are what
we shall use in general, but more important than these is the
relationship between the first and last expressions. We would
like it if the last one was similar to the first, but in some
sense simpler, so that we could generalize both statements to
one that can be proved inductively.

Certain similarities are immediately clear, as is the fact
that the last expression, if we fix $x$ and $x'$ rather than
taking the first expectation, involves functions of two variables 
rather than three, and a fourth power instead of an eighth
power. The only small difference is that we now have the 
function $G$ appearing rather than some arbitrary function
supported in $G$. This we shall have to incorporate into
our inductive hypothesis somehow.

However, in this small case, we can simply try to repeat the argument,
so let us continue with the calculation:
$$\eqalign{&\Bigl(\E_{y,z,t}f_{x,x'}(y,z)u_{x,x'}(y,t)
G(z,t)\Bigr)^4\cr
&=\Bigl(\E_{z,t}\E_yf_{x,x'}(y,z)u_{x,x'}(y,t)
G(z,t)\Bigr)^4\cr
&=\Bigl(\E_{z,t}\E_yf_{x,x'}(y,z)u_{x,x'}(y,t)
G_{x,x'}(z)G_{x,x'}(t)G(z,t)\Bigr)^4\cr
&\leq\Bigl(\E_{z,t}G_{x,x'}(z)G_{x,x'}(t)G(z,t)\Bigr)^2\Bigl(\E_{z,t}
\Bigl(\E_yf_{x,x'}(y,z)u_{x,x'}(y,t)G(z,t)\Bigr)^2\Bigr)^2\ .\cr}$$
Here, we used the fact that $f_{x,x'}(y,z)$ is non-zero only if
$G(x,z)$ and $G(x',z)$ are both equal to 1, with a similar statement
for $u_{x,x'}(y,t)$. We then applied the Cauchy-Schwarz inequality
together with the fact that $G$ squares to itself. Given that $G$
could be quite sparse, it was important here that we exploited its
sparseness to the full: with a lazier use of the Cauchy-Schwarz
inequality we would not have obtained the factor in the first bracket,
which will in general be small and not something we can afford to
forget about.

Now let us continue to manipulate the second bracket in the standard
way: expanding the inner square, rearranging, and applying Cauchy-Schwarz.
This time, in order not to throw away any sparseness information, we will 
bear in mind that the expectation over $y$ and $y'$ below is zero unless
all of $G(x,y)$, $G(x',y)$, $G(x,y')$ and $G(x',y')$ are equal to 1.
$$\eqalign{&\Bigl(\E_{z,t}\Bigl(\E_yf_{x,x'}(y,z)u_{x,x'}(y,t)
G(z,t)\Bigr)^2\Bigr)^2\cr
&=\Bigl(\E_{y,y'}G_{x,x',y,y'}\E_{z,t}f_{x,x',y,y'}(z)u_{x,x',y,y'}(t)G(z,t)\Bigr)^2\cr
&\leq\Bigl(\E_{y,y'}G_{x,x',y,y'}\Bigr)
\Bigl(\E_{y,y'}\Bigl(\E_{z,t}f_{x,x',y,y'}(z)u_{x,x',y,y'}(t)G(z,t)\Bigr)^2
\Bigr)\ .\cr}$$

We have now got down to functions of one variable, apart from the term
$G(z,t)$. Instead of worrying about this, let us continue the process.
$$\eqalign{\Bigl(\E_{z,t}f_{x,x',y,y'}(z)&u_{x,x',y,y'}(t)G(z,t)\Bigr)^2\cr
&=\Bigl(\E_t\E_zf_{x,x',y,y'}(z)u_{x,x',y,y'}(t)G(z,t)\Bigr)^2\ .\cr}$$

Now we shall apply Cauchy-Schwarz again, and again we must be careful
to use the full strength of the inequality by taking account that for
most values of $t$ the expectation over $z$ is zero. We can do this
by noting that 
$$u_{x,x',y,y'}(t)=u_{x,x',y,y'}(t)G_{x,x'}(t)G_{y,y'}(t)$$
so the last expression above is at most
$$\Bigl(\E_tG_{x,x'}(t)G_{y,y'}(t)\Bigr)
\Bigl(\E_t\Bigl(\E_zf_{x,x',y,y'}(z)u_{x,x',y,y'}(t)G(z,t)\Bigr)^2\Bigr).$$

The second term in this product is at most
$$\E_t\Bigl(\E_zf_{x,x',y,y'}(z)G_{x,x'}(t)G_{y,y'}(t)G(z,t)\Bigr)^2,$$
which equals
$$\E_t\E_{z,z'}f_{x,x',y,y',z,z'}G_{x,x'}(t)G_{y,y'}(t)G_{z,z'}(t).$$

Let us put all this together and see what the upper bound is that
we have obtained. It works out to be 
$$\eqalign{\Bigl(&\E_{y,z,t}G(y,z)G(y,t)G(z,t)\Bigr)^4\E_{x,x'}
\Bigl(\E_{z,t}G_{x,x'}(z)G_{x,x'}(t)G(z,t)\Bigr)^2
\Bigl(\E_{y,y'}G_{x,x',y,y'}\Bigr)\cr
&\E_{y,y'}\Bigl(\E_tG_{x,x'}(t)G_{y,y'}(t)\Bigr)\E_{z,z'}f_{x,x',y,y',z,z'}
\E_tG_{x,x'}(t)G_{y,y'}(t)G_{z,z'}(t).\cr}$$
Here we have been somewhat sloppy with our notation: a more correct
way of writing the above expression would be to have different names
for the variables in different expectations. If one does that and
then expands out the powers of the brackets, then one obtains an
expression with several further variables besides $x,x',y,y',z,z'$
and $t$. One takes the average, over all these variables, of an
expression that includes $f_{x,x',y,y',z,z'}$ and many terms 
involving the function $G$ applied to various pairs of the variables.
Recall that this is what we were trying to do.

We can interpret this complicated expression as follows. We allow the
variables to represent the vertices of a quadripartite graph $\Gamma$,
with two variables $q$ and $r$ joined by an edge if $G(q,r)$ appears
in the product. For example, the $G_{z,z'}(t)$ that appears at the end
of the expression is short for $G(z,t)G(z',t)$, so it would tell us
that $zt$ and $z't$ were edges of the graph (assuming that those
particular variables had not had their names changed).

When we assign values in $X$, $Y$, $Z$ and $T$ to the various
variables, we are defining a quadripartite map from the vertex
set of $\Gamma$ to the set $X\cup Y\cup Z\cup T$. And the
product of all the terms involving $G$ is telling us whether a particular
assignment to the variables of values in $X$, $Y$, $Z$ and $T$
results in a graph homomorphism from $\Gamma$ to $G$. 

Thus, the expression we obtain is an expectation over all such 
quadripartite maps $\phi$ of $f_{x,x',y,y',z,z'}$ multiplied by
the characteristic function of the event ``$\phi$ is a homomorphism.''

Notice that in this expression the function $f$ appears eight times,
as it does in the expression with which we started, since that
contains a single $f$ inside the bracket, which is raised to the
eighth power. This is important, as we need our inequality to scale
up in the right way. But equally important is that this scaling 
should occur correctly in $G$ as well. We can think of $G$ as put
together out of six functions (one for each pair of vertex sets).
Let us now reflect this in our notation, writing $G_{XY}$ for the
part of $G$ that joins $X$ to $Y$, and so on. If we want to make 
explicit the fact that $f$, $u$ and $w$ are zero except at triangles
in $G$, then we can rewrite the first expression as
$$\eqalign{\Bigl(\E_{x,y,z,t}f(x,y,z)u(x,y,t)w(y,z,t)G_{XY}(x,y)&G_{XZ}(x,z)
G_{XT}(x,t)\cr
&G_{YZ}(y,z)G_{YT}(y,t)G_{ZT}(z,t)\Bigr)^8.\cr}$$
This makes it clear that each part of $G$ (such as $G_{XY}$) occurs
eight times.  In order to have a useful inequality we need the
same to be true for the final expression that we are using to
bound this one. As it is written at the moment, $G_{XT}$,
$G_{YT}$ and $G_{ZT}$ are used eight times each, but $G_{XY}$,
$G_{YZ}$ and $G_{XZ}$ are used only four times each. However,
there are once again some implicit appearances, hidden in our
assumptions about when $f$ can be non-zero. In particular, we can afford
to multiply $f_{x,x',y,y',z,z'}$ by the product over all graph
terms, such as $G_{YZ}(y',z)$, that must equal 1 if 
$f_{x,x',y,y',z,z'}$ is non-zero. This gives us four extra 
occurrences of each of $G_{XY}$, $G_{YZ}$ and $G_{XZ}$.

We eventually want to show that if $\oct(f)$ is small and all the
functions such as $G_{XY}$ are ``sufficiently quasirandom'', then the
expression with which we started is small. In order to see what we do
next, let us abandon our current example, since it has become quite
complicated, and instead look at a simpler example that has the same
important features. In order to make this simpler example properly
illustrative of the general case, it will help if we no longer assume
that $G$ uses all the vertices in $X$, $Y$, $Z$ and $T$. Rather, we
shall let $P$, $Q$, $R$ and $S$ be subsets of $X$, $Y$, $Z$ and $T$,
respectively, and $G$ will be a graph that does not join any vertices
outside these subsets.  Then we shall consider how to approximate the
quantity
$$\E_{x,y,z,t}f(x,y,z)G(x,t)G(y,t)G(z,t)P(x)Q(y)R(z)S(t)$$
by the quantity
$$\E_{x,y,z,t}f(x,y,z)\d_{XT}G(y,t)G(z,t)P(x)Q(y)R(z)S(t),$$
where $\d_{XT}$ is now the {\it relative} density of $G$ inside
the set $P\times S$ (rather than its absolute density inside
$X\times T$). The sets $P$, $Q$, $R$ and $S$ will themselves
have densities, which we shall call $\d_X$, $\d_Y$, $\d_Z$
and $\d_T$.

To begin with, we define a function $g$ in the variables $x$
and $t$ by taking $g(x,t)$ to be $G(x,t)-\d_{XT}$ when 
$(x,t)\in P\times S$ and $0$ otherwise. The idea behind this
definition is that we want to subtract from $G(x,t)$ a function
that is supported in $P\times S$ and constant there, in such 
a way that the average becomes zero. Once we have done that,
our task is then to show that
$$\E_{x,y,z,t}f(x,y,z)g(x,t)G(y,t)G(z,t)P(x)Q(y)R(z)S(t)$$
is small, provided that $\oct(g)=\E_{x,x',t,t'}g_{x,x',t,t'}$ 
is small enough.

The technique of proof is the same as we have already seen: we give
the argument mainly to illustrate what we can afford to ignore and
what we must be careful to take account of. Since $g$ is a function
of two variables, we shall start with the expression
$$\eqalign{\Bigl(\E_{x,y,z,t}
&f(x,y,z)g(x,t)G(y,t)G(z,t)\Bigr)^4\cr
&=\Bigl(\E_{y,z,t}\E_xg(x,t)f(x,y,z)G(y,t)G(z,t)\Bigr)^4\cr
&\le\Bigl(\E_{y,z,t}G(y,z)G(y,t)G(z,t)\Bigr)^2
\Bigl(\E_{y,z,t}\Bigl(\E_xg(x,t)f(x,y,z)G(y,t)G(z,t)\Bigr)^2\Bigr)^2.\cr}$$
Now, we shall eventually be assuming that $\oct(g)$ is significantly 
smaller than the densities of any of the parts of $G$, but not necessarily
smaller than the densities of the sets $P$, $Q$, $R$ and $S$. The effect
on our calculations is that we can afford to throw away the $G$-densities
(by replacing them by 1) but must be careful to keep account of the
densities of vertex sets. Thus, we may replace the expectation
$\E_{y,z,t}G(y,z)G(y,t)G(z,t)$ in the first bracket by the larger 
expectation $\E_{y,z,t}Q(y)R(z)S(t)$. (This is of course easily seen to be
$\d_Y\d_Z\d_T$, but in more general situations it will not necessarily
be easy to calculate.) 

As for the second part of the product, it equals
$$\Bigl(\E_{y,z,t}G(y,t)G(z,t)\Bigl(\E_xg(x,t)f(x,y,z)\Bigr)^2\Bigr)^2,$$
which we can afford to bound above by
$$\eqalign{\Bigl(\E_{y,z,t}
&Q(y)R(z)S(t)\Bigl(\E_xg(x,t)f(x,y,z)\Bigr)^2\Bigr)^2\cr
&=\Bigl(\E_{x,x'}\E_{y,z,t}g_{x,x'}(t)f_{x,x'}(y,z)Q(y)R(z)S(t)\Bigr)^2\cr
&=\Bigl(\E_{x,x'}\E_{y,z,t}g_{x,x'}(t)P(x)P(x')f_{x,x'}(y,z)\Bigr)^2\cr
&\le\Bigl(\E_{x,x'}P(x)P(x')\Bigr)\Bigl(\E_{x,x'}\Bigl(\E_{y,z,t}
g_{x,x'}(t)f_{x,x'}(y,z)\Bigr)^2\Bigr).\cr}$$
Now we concentrate our efforts on the second bracket.
$$\eqalign{\Bigl(\E_{y,z,t}g_{x,x'}(t)&f_{x,x'}(y,z)\Bigr)^2\cr
&=\Bigl(\E_{y,z}Q(y)R(z)f_{x,x'}(y,z)\E_tg_{x,x'}(t)\Bigr)^2\cr
&\le\Bigl(\E_{y,z}Q(y)R(z)f_{x,x'}(y,z)^2\Bigr)
\Bigl(\E_{y,z}Q(y)R(z)\Bigl(\E_tg_{x,x'}(t)\Bigr)^2\Bigr).\cr}$$
Since $f$ is a function of three variables, we are even more prepared
to bound $f_{x,x'}(y,z)^2$ above by 1 than we were with $G$. That is,
we can bound the first bracket above by $\E_{y,z}P(x)P(x')Q(y)R(z)$.
The second equals $\E_{y,z,t,t'}Q(y)R(z)g_{x,x',t,t'}$. Since the
second is automatically zero if $P(x)P(x')$ is zero, we can even
afford to bound the first one by $\E_{y,z}Q(y)R(z)$.

Putting all this together, we find that
$$\Bigl(\E_{x,y,z,t}f(x,y,z)g(x,t)G(y,t)G(z,t)\Bigr)^4$$
is at most
$$\eqalign{\Bigl(\E_{y,z,t}Q(y)R(z)S(t)\Bigr)^2
&\Bigl(\E_{x,x'}P(x)P(x')\Bigr)\cr
&\Bigl(\E_{x,x'}\Bigl(\E_{y,z}Q(y)R(z)\Bigr)
\Bigl(\E_{y,z,t,t'}Q(y)R(z)g_{x,x',t,t'}\Bigr)\Bigr).\cr}$$
It is not hard to check that this equals $\d_X^2\d_Y^4\d_Z^4\d_T^2\oct(g)$.
This quantity will count as a small error if $\oct(g)$ is small
compared with $\d_X^2\d_T^2$, since then our upper bound is small
compared with its trivial maximum of $\d_X^4\d_Y^4\d_Z^4\d_T^4$
(which, in the general case, is rather less trivial).

An important point to note about the above argument is that even
though the expression we started with included a function of three
variables, it did not cause us any difficulty because we were
eventually able to bound it above in a simple way. This explains why
an inductive argument is possible: when we are dealing with functions
of $k$ variables $\seq x k$, we do not have any trouble from functions
of more variables, provided that at least one of $\seq x k$ is not
included in them. 

Of course, once we have replaced $G(x,t)$ by $\d_{XT}P(x)S(t)$ we can run 
similar arguments to replace $G(y,t)$ and $G(z,t)$ by $\d_{YT}Q(y)S(t)$
and $\d_{ZT}R(z)S(t)$, respectively. Thus, there will be three nested 
inductions going on at once: the number of variables $k$ in the function under 
consideration, the number of functions of $k$ variables still left to 
consider, and the number of steps taken in the process of replacing a
function $f$ by a function of the form $f_{x_1,x_1',\dots,x_k,x_k'}$.
Section 4 is concerned with the last of these, and the first two are
dealt with in Section 5. 

\bigskip

\noindent {\bf \S 3. Some basic definitions.}
\medskip

The need for a more compact notation should by now be clear. In
this section, we shall provide such a notation and also explain 
the terminology that will be needed to state our main results.
\bigskip

\noindent {3.1. {\sl Hypergraphs and chains.}}
\medskip

An {\it $r$-partite hypergraph} is a sequence $\seq X r$ of disjoint
sets, together with a collection $\ch$ of subsets $A$ of
$X_1\cup\dots\cup X_r$ with the property that $|A\cap X_i|\le 1$ for
every $i$. The sets $X_i$ are called {\it vertex sets} and their
elements are {\it vertices}. The elements of $\ch$ are called {\it
edges}, or sometimes {\it hyperedges} if there is a danger of
confusing them with edges in the graph-theoretic sense. A hypergraph
is $k$-{\it uniform} if all its edges have size $k$. (Thus, a
2-uniform hypergraph is a graph.)

An $r$-partite hypergraph $\ch$ is called an $r$-{\it partite chain}
if it has the additional property that $B$ is an edge of $\ch$
whenever $A$ is an edge of $\ch$ and $B\subset A$. Thus, an
$r$-partite chain is a particular kind of combinatorial simplicial
complex, or down-set. Our use of the word ``chain'' is non-standard
(in particular, it has nothing to do with the notion of a chain
complex in algebraic topology). We use it because it is quicker to
write than ``simplicial complex''.

If the largest size of any edge of $\ch$ is $k$, then we shall sometimes
say that $\ch$ is a $k$-{\it chain}. 

\medskip

\noindent {3.2. {\sl Homomorphisms and $r$-partite functions.}}
\medskip

Let $\seq E r$ and $\seq X r$ be two sequences of disjoint finite sets. 
If $\phi$ is a map from $E_1\cup\dots\cup E_r$ to $X_1\cup\dots\cup X_r$
such that $\phi(E_i)\subset X_i$ for every $i$, we shall say that $\phi$
is an $r$-{\it partite function}. 

Let $\cj$ be an $r$-partite chain with vertex sets $\seq E r$ and
let $\ch$ be an $r$-partite chain with vertex sets $\seq X r$. Let
$\phi$ be an $r$-partite function from the vertices of $\cj$ to the
vertices of $\ch$. We shall say that $\phi$ is a {\it homomorphism}
from $\cj$ to $\ch$ if $\phi(A)\in\ch$ whenever $A\in\cj$. We shall
write $\Hom(\cj,\ch)$ for the set of all homomorphisms from $\cj$
to~$\ch$.
\medskip
\noindent {3.3. {\sl $A$-functions and $\cj$-functions.}}
\medskip

Let $\Phi$ be the set of all $r$-partite maps from $E_1\cup\dots\cup E_r$ 
to $X_1\cup\dots\cup X_r$. We shall also consider some special classes of 
functions defined on $\Phi$. If $A$ is a subset of $E_1\cup\dots\cup E_r$ 
such that $|A\cap E_i|\le 1$ for every $i$, then a function 
$f:\Phi\ra[-1,1]$ will be called an $A$-{\it function} if the value
of $f(\phi)$ depends only on the image $\phi(A)$. If $\cj$ is an $r$-partite
chain with vertex sets $\seq E r$, then a $\cj$-{\it function} is a function
$f:\Phi\ra[-1,1]$ that can be written as a product $f=\prod_{A\in\cj}f^A$,
where each $f^A$ is an $A$-function.

The definition of $A$-functions and $\cj$-functions is introduced in
order to deal with situations where we have a function of several variables 
that can be written as a product of other functions each of which depends
on only some of those variables. We met various functions of this type
in the previous section. Let us clarify the definition with another small 
example. Suppose that we have three sets $X_1$, $X_2$ and $X_3$ and a function 
$f:X_1^2\times X_2\times X_3\ra[-1,1]$ of the form
$$f(x_1,x_1',x_2,x_3)=f_1(x_1,x_2)f_2(x_1,x_3)f_3(x_1',x_2)f_4(x_1',x_3)\ .$$
Let $E_1=\{1,1'\}$, $E_2=\{2\}$ and $E_3=\{3\}$. There is an obvious 
one-to-one correspondence between quadruples $(x_1,x_1',x_2,x_3)$ and 
tripartite maps from $E_1\cup E_2\cup E_3$: given such a sequence one
associates with it the map $\phi$ that takes $1$ to $x_1$, $1'$ to $x_1'$,
$2$ to $x_2$ and $3$ to $x_3$. Therefore, we can if we wish change to a 
more opaque notation and write
$$f(\phi)=f_1(\phi)f_2(\phi)f_3(\phi)f_4(\phi)\ .$$
Now $f_2(\phi)=f_2(\phi(1),\phi(3))=f_2\big(\phi(\{1,3\})\bigr)$, so
$f_2$ is a $\{1,3\}$-function. Similar remarks can be made about
$f_1$, $f_3$ and $f_4$. It follows that $f$ is a $\cj$-function if
we take $\cj$ to be the chain consisting of the sets $\{1,2\}$,
$\{1,3\}$, $\{1',2\}$ and $\{1',3\}$ and all their subsets. The
fact that the subsets are not mentioned in the formula does not
matter, since if $C$ is one of these subsets we can take the function 
that is identically 1 as our $C$-function.

An important and more general example is the following. As above, let
$\cj$ be an $r$-partite chain with vertex sets $\seq E r$ and let
$\ch$ be an $r$-partite chain with vertex sets $\seq X r$. For each
$\phi$ in $\Phi$ and each $A\in\cj$ let $H^A(\phi)$ equal 1 if
$\phi(A)\in\ch$ and $0$ otherwise. Let
$\ch(\phi)=\prod_{A\in\cj}H^A(\phi)$.  Then $\ch(\phi)$ equals $1$ if
$\phi\in\Hom(\cj,\ch)$ and $0$ otherwise.  In other words, the
characteristic function of $\Hom(\cj,\ch)$ is a $\cj$-function. We
stress that $\ch(\phi)$ depends on $\cj$; however, it is convenient to
suppress this dependence in the notation.  Our counting lemma will
count homomorphisms from small chains $\cj$ to large quasirandom
chains $\ch$, so we can regard our main aim as being to estimate the
sum (or equivalently, expectation) of $\ch(\phi)$ over all
$\phi\in\Phi$. However, in order to do so we need to consider more
general $\cj$-functions.

The $\cj$-functions we consider will be supported in a chain $\ch$ in
the following sense. Let us say that an $A$-function $f^A$ is {\it
supported in} $\ch$ if $f^A(\phi)$ is zero whenever $\phi(A)$ fails to
be an edge of $\ch$. Equivalently, $f^A$ is supported in $\ch$ if
$f^A=f^AH^A$, where $H^A$ is as defined above.  We shall say that $f$
is a $\cj$-{\it function on} $\ch$ if it can be written as a product
$\prod_{A\in\cj}f^A$, where each $f^A$ is an $A$-function supported in
$\ch$. If $f$ is a $\cj$-function on $\ch$, then $f(\phi)=0$ whenever
$\phi$ does not belong to $\Hom(\cj,\ch)$. That is,
$f(\phi)=f(\phi)\ch(\phi)$. Notice that the product of any $\cj$
function with the function $\ch$ will be a $\cj$-function on $\ch$.

This is another definition that came up in the previous section.  In
that case, the three functions in the product
$f(x,y,z)u(x,y,t)v(y,z,t)$ considered in the previous section were all
supported in the chain $\ch$ that consisted of the triangles in the
graph $G$, the edges of $G$, and the vertices of $G$. If we let $\cj$
be the chain consisting of the sets $\{x,y,z\}$, $\{x,y,t\}$, $\{y,z,t\}$
and all their subsets (where we are regarding the letters as names of
variables rather than elements of $X$, $Y$, $Z$ and $T$), then this
product is a $\cj$-function on $\ch$. 
\medskip

\noindent {3.4. {\sl The index of a set, and relative density in a chain.}}
\medskip

Let $\ch$ be an $r$-partite chain with vertex sets $\seq X r$. Given a
set $F\in\ch$, define its {\it index} $i(F)$ to be the set of all $i$
such that $F\cap X_i$ is non-empty. (Recall that $F\cap X_i$ is a
singleton for each such $i$.) For any set $A$ in any $r$-partite
chain, let $H(A)$ be the collection of all sets $E\in\ch$ of index
equal to that of $A$. If $A$ has cardinality $k$, then let $H_*(A)$ be
the collection of all sets $D$ of index $i(A)$ such that $C\in\ch$
whenever $C\subset D$ and $C$ has cardinality $k-1$. (Since $\ch$ is a
chain, it follows from this that all proper subsets of $D$ belong to
$\ch$. Note that we do not require $D$ to belong to $\ch$.) Clearly
$H(A)\subset H_*(A)$. The {\it relative density of} $H(A)$ {\it in}
$\ch$ is defined to be $|H(A)|/|H_*(A)|$. We will denote it by $\d_A$.

Once again, the example in the last section illustrates the importance
of $H_*(A)$. Let us rename the vertex sets $X$, $Y$, $Z$ and $T$ as
$X_1,X_2,X_3$ and $X_4$. If $\ch$ is a 3-chain that consists of the
edges and vertices of the graph $G$, and some collection of triangles
of $G$, and if $A=\{1,2,3\}$, say, then $H_*(A)$ consists of all
triangles in $G$ with one vertex in each of $X_1$, $X_2$ and $X_3$,
while $H(A)$ consists of all 3-edges of $\ch$ with one vertex in each
of $X_1$, $X_2$ and $X_3$. Thus, $\d_A$ measures the proportion of the
triangles in $G$ that are edges in $\ch$.

It is useful to interpret the relative density $\d_A$
probabilistically: it is the conditional probability that a randomly
chosen set $D\subset X_1\cup\dots\cup X_r$ of index $i(A)$ belongs to
$\ch$ (and hence to $H(A)$), given that all its proper subsets belong
to $\ch$.
\medskip

\noindent {\bf Notational remark.} It may help the reader to remember
the definitions in this section if we explicitly point out that most
of the time we are adopting the following conventions. The symbols
$\cj$ and $\ck$ are used for chains of fixed size that are embedded
into a chain $\ch$ of size tending to infinity. From these we
sometimes form other chains: for instance, $\cj_1$ will be a chain of
fixed size derived from a chain $\cj$, and $\ch(x)$ will be a
chain of size tending to infinity that depends on a point $x$. The
letter $H$ will tend to be reserved for set systems connected with
$\ch$ where the sets all have the same index. The same goes for
functions derived from $\ch$. For example, we write $\ch(\phi)$
because we use the full chain $\ch$ to define the function, whereas we
write $H^A(\phi)$ because for that we just use sets of index $i(A)$,
which all have size $|A|$. Similarly, we write $H_*(A)$ because all 
sets in $H_*(A)$ have index $i(A)$.
\medskip

\noindent {3.5. {\sl $\oct(f^A)$ for an $A$-function $f^A$.}}
\medskip
We are building up to a definition of quasirandomness for $\ch(A)$.
An important ingredient of the definition is a weighted count of
combinatorial octahedra, which generalizes the definition introduced
in the last section. If $f$ is a function of three variables $x$, $y$
and $z$ that range over sets $X$, $Y$ and $Z$, respectively, then we
defined $\oct(f)$ to be $\E_{x,x',y,y',z,z'}f_{x,x',y,y',z,z'}$. In 
full, this is the expectation over all $x,x'\in X$,
$y,y'\in Y$ and $z,z'\in Z$ of
$$f(x,y,z)f(x,y,z')f(x,y',z)f(x,y',z')
f(x',y,z)f(x',y,z')f(x',y',z)f(x',y',z')\ .$$
Similarly, if $f$ is a function of $k$ variables $\seq x k$,
with each $x_i$ taken from a set $X_i$, then 
$$\oct(f)=\E_{x_1^0,x_1^1\in X_1}\dots\E_{x_k^0,x_k^1\in X_k}
\prod_{\e\in\{0,1\}^k}f(x_1^{\e_1},\dots,x_k^{\e_k})\ .$$
In the spirit of the previous section, we can (and shall) also
write this as $\E_\sigma f_\sigma$, where $\sigma$ is shorthand
for $x_1,x_1',\dots,x_k,x_k'$.  

To give a formal definition in more general situations it is
convenient to use the language of $A$-functions, though in fact
we shall try to avoid this by assuming without loss of generality
that the set $A$ we are talking about is the set $\{1,2,\dots,k\}$. 
Nevertheless, here is the definition. As before, let $\cj$
and $\ch$ be $r$-partite chains with vertex sets $\seq E r$ and $\seq
X r$, let $\Phi$ be the set of all $r$-partite maps from
$E_1\cup\dots\cup E_r$ to $X_1\cup\dots\cup X_r$ and let $A\in\cj$. We
can think of an $A$-function as a function defined on the product of
those $X_i$ for which $i\in i(A)$.  However, we can also think of it
as a function $f^A$ defined on $\Phi$ such that $f^A(\phi)$ depends
only on $\phi(A)$. To define $\oct(f^A)$ in these terms, we construct
a set system $\cb$ as follows. Let $k$ be the cardinality of the set
$A$. For each $i\in i(A)$ let $U_i$ be a set of cardinality 2, let $U$
be the union of the $U_i$ (which we suppose to be disjoint) and let
$\cb$ consist of the $2^k$ sets $B\subset U$ such that $|B\cap U_i|=1$
for every $i$. Let $\Omega$ be the set of all $k$-partite maps
$\omega$ from $\bigcup_{i\in i(A)}U_i$ to $\bigcup_{i\in i(A)}X_i$
(meaning that $\omega(U_i)\subset X_i$ for every $i\in i(A)$).

We now want to use $f^A$, which is defined on $\Phi$, to define a
$B$-function $f^B$ on $\Omega$, for each $B\in\cb$.  There is only one
natural way to do this. Given $\omega\in\Omega$ and $B\in\cb$, we
would like $f^B(\omega)$ to depend on $\omega(B)$; we know that $B$
and $\omega(B)$ have the same index as $A$; so we choose some
$\phi\in\Phi$ such that $\phi(A)=\omega(B)$ and define $f^B(\omega)$
to be $f^A(\phi)$. This is well-defined, since if $\phi(A)=\phi'(A)$,
then $f^A(\phi)=f^A(\phi')$, because $f^A$ is an $A$-function.

We now define
$$\oct(f^A)=\E_{\omega\in\Omega}\prod_{B\in\cb}f^B(\omega)\ .$$
Let us see why this agrees with our earlier definition. There,
for simplicity, we took $A$ to be the set $\{1,2,\dots,k\}$. 
Then for each $i\le k$ we let $U_i=\{x_i^0,x_i^1\}$, and $\cb$
consisted of all sets of the form $B_\e=\{x_1^{\e_1},\dots,x_k^{\e_k}\}$,
with $\e=(\e_1,\dots,\e_k)\in\{0,1\}^k$. The set $\Omega$ was the
set of all ways of choosing $x_i^0$ and $x_i^1$ in $X_i$, for
each $i\le k$. (Again there is a deliberate ambiguity in our
notation. When we say that $U_i=\{x_i^0,x_i^1\}$ we are thinking
of $x_i^0$ and $x_i^1$ as symbols for variables, and when we
choose elements of $X_i$ with those names, we are thinking of
this choice as a function from the set $\{x_i^0,x_i^1\}$ of symbols
to the set $X_i$.) Given $\omega\in\Omega$ and $B=B_\e\in\cb$, we
have to define $f^{B_\e}(\omega)$. In principle a function of $\omega$
can depend on all the variables $x_i^0$ and $x_i^1$, but $f^{B_\e}$ is
a $B_\e$-function, and therefore depends just on the variables
$x_i^{\e_i}$. Now $\Phi$ can be thought of as the set of ways of
choosing $y_i\in X_i$ for each $i\le k$. In other words, we regard
$A$ as the set of variables $\{y_1,\dots,y_k\}$ and $\phi$ as a way
of assigning values to these variables. Thus, to define $f^{B_\e}(\omega)$
we choose $\phi$ such that $\phi(A)=\omega(B_\e)$, which means
that $\phi(y_i)$ must equal $\omega(x_i^{\e_i})$ for each $i$. 
(Equivalently, thinking of $y_i$ and $x_i^{\e_i}$ as the assigned
values, it means merely that $x_i^{\e_i}$ must equal $y_i$.) 
But then $f(\phi)=f(y_1,\dots,y_k)=f(x_1^{\e_1},\dots,x_k^{\e_k})$. 
And now it is clear that the two expressions for $\oct(f)$ denote
the same quantity.
\medskip
\noindent {3.6. {\sl Octahedral quasirandomness.}}
\medskip

We come now to the first of two definitions that are of great
importance for this paper. Let $\ch$ be a chain, let $f^A$ be an
$A$-function, for some $A$ that does not necessarily belong to $\ch$,
and suppose that $f^A$ is supported in $H_*(A)$, in the sense that
$f^A(\phi)=0$ whenever $\phi(A)\notin H_*(A)$.  Equivalently, suppose
that whenever $f^A(\phi)\ne 0$ we have $\phi(C)\in\ch$ for every
proper subset $C\subset A$. Loosely speaking, we shall say that $f$ is
{\it octahedrally quasirandom relative to} $\ch$ if $\oct(f^A)$ is
significantly smaller than one might expect.

To turn this idea into a precise definition, we need to decide what we
expect. Let $\cb$ be the set system defined in the previous
subsection. If $B\in\cb$, then $f^B(\omega)$ is defined to be the
value of $f^A(\phi)$ for any $\phi$ with $\phi(A)=\omega(B)$. If
$f^B(\omega)\ne 0$, then $f^A(\phi)\ne 0$ so $\phi(A)\in H_*(A)$, by
assumption, and hence $\omega(B)\in H_*(A)$. Therefore, a necessary
condition for $\prod_{B\in\cb}f^B(\omega)$ to be non-zero is that
$\omega(D)\in\ch$ for every $D$ that is a proper subset of some
$B\in\cb$. Let $\ck'$ be the chain consisting of all such sets. Thus,
$\ck'$ consists of all subsets of $U_1\cup\dots\cup U_k$ that intersect
each $U_i$ in at most a singleton and do not intersect every $U_i$. Then,
since $|f^B(\omega)|\le 1$ for every $B$ and every $\omega$, a trivial
upper bound for $\oct(f^A)$ is
$$\E_{\omega\in\Omega}\prod_{D\in\ck'}H^D(\omega)\ ,$$
which we shall call $\oct(H_*(A))$, since it counts the number of
(labelled, possibly degenerate) combinatorial $k$-dimensional
octahedra in $H_*(A)$. 

We could if we wanted declare $\oct(f^A)$ to be small if it is small
compared with $\oct(H_*(A))$. Instead, however, since we shall be
working exclusively with quasirandom chains, it turns out to be more
convenient to work out how many octahedra we expect $H(A)$ to have,
given the various relative densities, and use that quantity for
comparison. (It might seem more natural to use $H_*(A)$, but for
the particular functions $f^A$ that we shall need to consider, 
$\oct(f^A)$ will tend to be controlled by the smaller quantity
$\oct(H(A))$. But in the end this is not too important because
when we are looking at $\oct(f^A)$ we think of the density 
$\d_A$ as ``large''.) 

Let us therefore write $\ck$ for the set
of {\it all} subsets of sets in $\cb$ (so $\ck=\cb\cup\ck'$).
It is helpful to recall the interpretation of relative
densities as conditional probabilities. Suppose that we choose
$\omega$ randomly from $\Omega$, and also that $\ch$ behaves in 
a random way. Then the probability that $H^D(\omega)=1$ given that
$H^C(\omega)=1$ for every $C\subsetneq D$ is the probability that
$\omega(D)\in\ch$ given that $\omega(C)\in\ch$ for every $C\subsetneq D$,
which is $\d_D$. Because $\ch$ behaves randomly, we expect all these
conditional probabilities to be independent, so we expect that 
$\E_{\omega\in\Omega}\prod_{D\in\ck}H^D(\omega)$ will be approximately
$\prod_{D\in\ck}\d_D$. Accordingly, we shall say that $f^A$ is
$\eta$-{\it octahedrally quasirandom} if
$$\oct(f^A)\le\eta\prod_{D\in\ck}\d_D\ .$$
Since octahedral quasirandomness is the only form of quasirandomness
that we use in this paper, we shall often omit the word ``octahedrally''
from this definition.

It is not necessary to do so, but one can rewrite the right-hand 
side more explicitly. For each subset $C\subset A$, 
there are $2^{|C|}$ sets $D\in\ck$
with the same index as $C$. (We can think of these as $|C|$-dimensional
faces of the octahedron with index $i(C)$.) Therefore, 
$$\eta\prod_{D\in\ck}\d_D
=\eta\prod_{C\subset A}\d_C^{2^{|C|}}\ .$$

The main use of the definition of quasirandomness for $A$-functions
is to give us a precise way of saying what it means for a $k$-partite
$k$-uniform hypergraph to ``sit quasirandomly inside a $k$-partite
$(k-1)$-chain''. Let $A$ and $\ch$ be as above. The $k$-uniform
hypergraph we would like to discuss is $H(A)$. Associated with 
this hypergraph is its ``characteristic function'' $H^A$ and its 
relative density $\d_A$. The $(k-1)$-chain is the set of all edges
of $\ch$ with index some proper subset of $A$. Define an $A$-function $f^A$
by setting $f^A(\phi)$ to equal $H^A(\phi)-\d_A$ if $\phi(A)\in H_*(A)$
and zero otherwise. An important fact about $f^A$ is that its average
is zero. To see this, note that $f^A(\phi)=H(\phi(A))-\d_A$ when 
$\phi(A)\in H_*(A)$ and $f^A(\phi)=0$ otherwise. Therefore, the average
over all $\phi$ such that $\phi(A)\notin H_*(A)$ is trivially zero, 
while the average over all $\phi$ such that $\phi(A)\in H_*(A)$ is
zero because $\d_A$ is the relative density of $H(A)$ in $H_*(A)$.

We shall say that $H(A)$ is $\eta$-{\it octahedrally quasirandom}, or
just $\eta$-{\it quasirandom}, relative to $\ch$, if the function $f^A$
is $\eta$-quasirandom according to the definition given earlier.
The counting lemma, which we shall prove in \S 5, will show that
if $\ch$ is an $r$-partite chain and all its different parts of the
form $H(A)$ are quasirandom in this sense, then $\ch$ behaves like
a random chain with the same relative densities.

\medskip
\noindent {3.7. {\sl Quasirandom chains.}}
\medskip

We are now ready for the main definition in terms of which our
counting and regularity lemmas will be stated. Roughly speaking, a
chain $\ch$ is quasirandom if $H(A)$ is highly quasirandom relative
to $\ch$. However, there is an important subtlety to the definition,
which is that when we apply it we do so in situations where the
relative densities $\d_A$ tend to be very much smaller when the
sets $A$ are smaller, as we saw in the second example of the 
previous section. For this reason, we need to make much stronger
quasirandomness assumptions about $H(A)$ when $A$ is small, and it
is also very important which of these assumptions depend on which
densities. The full details of the following definition are not
too important -- they are chosen to make the proof work -- but the
dependences certainly are. 

One other comment is that our definition depends on a chain $\cj$. 
This is useful for an inductive hypothesis later. Roughly, if $\ch$
is quasirandom with respect to $\cj$ then $\cj$ embeds into $\ch$ in 
the expected way. Thus, the bigger $\cj$ is, the stronger the
statement. 

Now let us turn to the precise definition. Suppose that $\cj$ and
$\ch$ are $r$-partite chains. For each $A\in\cj$, let the relative
density of $H(A)$ in $\ch$ be $\d_A$ and suppose that $H(A)$ is
relatively $\eta_A$-quasirandom. Define a sequence 
$\e_k,\e_{k-1},\dots,\e_1$ by taking $\e_k=\e$ and
$$\e_{k-j}=2^{-jk-1}|\cj|^{-1}\Bigl(\e_{k-j+1}
\prod_{A\in\cj\atop|A|\ge k-j+1}\d_A\Bigr)^{2^{jk}}$$
when $j\ge 1$. Let $\eta_{k-j}$ be defined by the formula
$$\eta_{k-j}=(1/2)
\Bigl(\e_{k-j}\prod_{A\in\cj\atop|A|\ge k-j}\d_A\Bigr)^{2^{k(j+1)}}$$
for each $j$. Then $\ch$ is $(\e,\cj,k)$-{\it quasirandom}
if, for every $A\in\cj$ of size $j\le k$, we have the inequality
$\eta_A\le\eta_j$, or in other words $H(A)$ is $\eta_j$-quasirandom
relative to $H_*(A)$.

The parameter $k$ is also there just for convenience in our eventual
inductive argument. The counting lemma will imply that if $\phi$ is a random 
$r$-partite map from $\cj$ to an $(\e,\cj,k)$-quasirandom chain $\ch$, 
and if all sets in $\cj$ have size at most $k$, then the probability
that $\phi$ is a homomorphism differs from $\prod_{A\in\cj}\d_A$ by
at most $\e|\cj|\prod_{A\in\cj}\d_A$.
\bigskip


\noindent {\bf \S 4. The main lemma from which all else follows.}
\medskip

Before we tackle our main lemma it will help to prepare for it in
advance with a small further discussion of terminology. Let $\ch$ be
an $r$-partite chain with vertex sets $X_1,\dots,X_r$. Let $t\ge r$
and let $\seq x t$ be variables such that $x_i$ ranges over $X_i$ when
$i\le r$ and over some other $X_j$ if $i>r$. For each $j\le r$ let
$E_j$ be the set of $i$ such that $x_i$ ranges over $X_j$ (so, in
particular, $i\in E_i$ when $i\le r$).

Now let $\cj$ be an $r$-partite chain with vertex sets 
$\seq E r$. Suppose that the set $\{1,2,\dots,k\}$ does not
belong to $\cj$ but that all its proper subsets do.

We shall write $\tau$ for the sequence $(\seq x t)$. Note
that there is a one-to-one correspondence between such sequences
and $r$-partite maps from $E_1\cup\dots\cup E_r$ to $X_1\cup\dots\cup X_r$,
so we can also think of $\tau$ as such a map.



Our aim will be to find an upper bound for the modulus of a quantity of 
the form
$$\E_{\tau}f(\tau)\prod_{A\in\cj}g^A(\tau),$$
where $f$ is any function from $X_1\times\dots\times X_r$ to $\R$, and 
each $g^A$ is an $A$-function supported in $\ch$ and taking values in
$[-1,1]$. By $f(\tau)$ we mean $f(x_1,\dots,x_r)$, but for convenience
we add in the other variables on which $f$ does not depend.

In order to shorten the statement of the next lemma, let us describe
in advance a chain $\ck$ that appears in its conclusion. For each $i\le t$
we shall have a set $W_i$ of the form $\{i\}\times U_i$, where $U_i$ 
is a finite subset of $\N$. The chain $\ck$ will be an $r$-partite chain
with vertex sets $\seq F r$, where $F_j=\bigcup_{i\in E_j}W_i$. We shall
use the vertices of $\ck$ to index variables as follows: the element
$(i,h)$ of $W_i$ indexes a variable that we shall call $x_i^h$. When $i\le k$
the sets $U_i$ will be chosen in such a way that $(i,0)$ and $(i,1)$
both belong to $U_i$: it will sometimes be convenient to use the
alternative names $x_i$ and $x_i'$ for $x_i^0$ and $x_i^1$.

We shall use the letter $\omega$ to stand for the sequence of all
variables $x_i^j$, enumerated somehow. Equivalently, we can think of
$\omega$ as an $r$-partite map from $F_1\cup\dots\cup F_r$ to
$X_1\cup\dots\cup F_r$.

Let $\sigma$ be shorthand for the sequence $x_1,x_1',x_2,x_2',\dots,x_k,x_k'$.
Generalizing the notation from \S 2, if $f:X_1\times\dots\times X_r\ra\R$
we shall write $f_\sigma(\omega)$ for the expression
$\prod_{\e\in\{0,1\}^k}f(x_1^{\e_1},\dots,x_k^{\e_k},x_{k+1},\dots,x_r)$.
Once again, $\omega$ contains many more variables than the ones that 
appear in this expression, but since $f$ does not depend on them the
notation is unambiguous. (In fact, when we come to apply the lemma, 
$f$ will not even depend on $x_{k+1},\dots,x_r$.)
\bigskip


\noindent {\bf Lemma {4.1}.} {\sl Let the chains $\ch$ and $\cj$ be 
as just described. Then there is a chain $\ck$ of the kind that has
also just been described, with the following properties.

(i) Every set in $\ck$ has cardinality less than $k$. 

(ii) Let $\g:F_1\cup\dots\cup F_r\ra E_1\cup\dots\cup E_r$ be
the $r$-partite map $(i,j)\mapsto i$. (That is, for each $i\le t$,
$\g$ takes the elements of $W_i$ to $i$.) Then $\g$ is a homomorphism
from $\ck$ to $\cj$, and for each $A\in\cj$
of cardinality less than $k$ there are precisely $2^k$ sets 
$B\in\ck$ such that $\g(B)=A$. 

(iii) If $f$ is any function from $X_1\times\dots\times X_r$ to $\R$
and each $g^A$ is an $A$-function supported in $\ch$ and taking values 
in $[-1,1]$, then we have the inequality
$$\Bigl(\E_{\tau}f(\tau)\prod_{A\in\cj}g^A(\tau)\Bigr)^{2^k}\leq
\E_\omega f_{\sigma}(\omega)\prod_{B\in\ck}H^B(\omega)\ .$$}
\bigskip





\Proof  We shall prove this result by induction. To do this we shall
show that for each $j\leq k$ the 
left-hand side can be bounded above by a quantity of the following
form, which we shall write first and then interpret:
$$\E_{\omega_j}\Bigl(\prod_{A\in\ck_j}H^A(\omega_j)\Bigr)
\Bigl(\E_{\tau_j}f_{\sigma_j}(\tau_j)
\prod_{[j]\subset A}(g^A)_{\sigma_j}(\tau_j)
\prod_{[j]\not\subset A\atop |A|<k}(H^A)_{\sigma_j}(\tau_j)\Bigr)^{2^{k-j}}.$$
The set system $\ck_j$ here is a chain. Each vertex of $\ck_j$ belongs
to a set $V_i^j$ of the form $\{i\}\times U_i^j$ for some $i\le t$ and some 
finite subset $U_i^j$ of $\N$. The vertices are partitioned into
$r$ sets $E_1^j,\dots,E_r^j$, where $E_i^j=\bigcup_{h\in E_i}V_h^j$.
As before, $x_h^q$ stands for a variable indexed by the pair $(h,q)\in V_h^j$.
In the back of our minds, we identify $(i,0)$ with $i$ when $i\le r$:
in particular, we shall sometimes write $x_i$ instead of $x_i^0$, and if
$j\le k$ we shall sometimes write $[j]$ for the set 
$\{(1,0),(2,0),\dots,(j,0)\}$ rather than the more usual $\{1,2,\dots,j\}$.
We shall also sometimes write $x_i'$ for $x_i^1$.

For the products in the second bracket we have not mentioned the
condition $A\in\cj$, which always applies. In other words, the
products are over all sets $A\in\cj$ that satisfy the conditions
specified underneath the product signs. We write $\sigma_j$ as
shorthand for $(x_1,x_1',\dots,x_j,x_j')$. We also write $\tau_j$ for
the sequence $(x_{j+1},\dots,x_t)$. We define the sets $V_i^j$ in such
a way that $V_i^0$ is the singleton $\{(i,0)\}$ and is a subset of
each $V_i^j$: it is only the first bracket that depends on the new
variables. Finally, $\omega_j$ is an enumeration of all the variables
that are not included in $\tau_j$.


We shall not specify what the edges of the chain $\ck_j$ are (though
in principle it would be possible to specify them exactly), since
all that concerns us is that the map $\g$ that takes $(i,0)$ to $i$
is a homomorphism from $\ck_j$ to $\cj$ such that,
for each $A\in\cj$ of cardinality less than $k$, the number of sets
$B\in\ck_j$ with $\g(B)=A$ is $2^k-2^{k-j+|A\cap[j]|}$ if
$A\not\subset[j]$ and $2^k-2^{|A|}$ if $A\subset[j]$.

Let us explain these last numbers. They are what we need for the
inequality to be properly homogeneous in the way that we discussed in
\S 2. To see why they are the correct numbers, let us think about a
function of the form
$(H^A)_{\sigma_j}=(H^A)_{x_1,x_1',\dots,x_j,x_j'}$.  For each $i\le j$
such that $i\notin A$, there is no dependence of
$(H^A)_{\sigma_j}(\tau_j)$ on $x_i$ or $x_i'$, so in order for
$(H^A)_{\sigma_j}(\tau_j)$ not to be zero, the number of distinct sets
that are required to belong to $\ch$ is $2^{|A\cap[j]|}$. When we
raise to the power $2^{k-j}$, this must happen $2^{k-j}$ times, all
independently, except that if $A\subset[j]$ then $H^A$ does not depend
on any of the variables in $\tau_j$ so it needs to happen just
once. Thus, the number of sets required to be in $\ch$ is
$2^{k-j}2^{|A\cap[j]|}=2^{k-j+|A\cap[j]|}$ when $A\not\subset[j]$, and
it is $2^{|A\cap[j]|}=2^{|A|}$ when $A\subset[j]$. This falls short of
$2^k$ and the difference must be made up for in the first bracket.


Now that we have discussed the inductive hypothesis in detail, let us prove
it by repeating once again the basic technique: isolate one variable
and sum over it last, apply Cauchy-Schwarz carefully, expand out
a square, rearrange, and apply Cauchy-Schwarz carefully again.

As we did repeatedly in \S 2, we shall leave the first bracket
and concentrate on the second. That is, we shall find an upper
bound for 
$$\Bigl(\E_{\tau_j}f_{\sigma_j}(\tau_j)
\prod_{[j]\subset A}(g^A)_{\sigma_j}(\tau_j)
\prod_{[j]\not\subset A\atop |A|<k}(H^A)_{\sigma_j}(\tau_j)\Bigr)^{2^{k-j}}.$$
Let us write $\tau_j$ as $(x_{j+1},\tau_{j+1})$. The quantity above equals
$$\Bigl(\Bigl(\E_{\tau_{j+1}}\E_{x_{j+1}}f_{\sigma_j}(x_{j+1},\tau_{j+1})
\prod_{[j]\subset A}(g^A)_{\sigma_j}(x_{j+1},\tau_{j+1})
\prod_{[j]\not\subset A\atop |A|<k}(H^A)_{\sigma_j}(x_{j+1},\tau_{j+1})
\Bigr)^2\Bigr)^{2^{k-j-1}}.$$
Applying Cauchy-Schwarz, we find that this is at most the product of
$$\Bigl(\E_{\tau_{j+1}}\!\prod_{[j]\subset A\atop j+1\notin A}
(g^A)_{\sigma_j}(x_{j+1},\tau_{j+1})^2
\!\prod_{[j]\not\subset A\atop{j+1\notin A\atop|A|<k}}(H^A)_{\sigma_j}
(x_{j+1},\tau_{j+1})\Bigr)^{2^{k-j-1}}$$
and
$$\Bigl(\E_{\tau_{j+1}}\Bigl(\E_{x_{j+1}}f_{\sigma_j}(x_{j+1},\tau_{j+1})
\!\prod_{[j+1]\subset A}(g^A)_{\sigma_j}(x_{j+1},\tau_{j+1})
\!\prod_{[j+1]\not\subset A\atop |A|<k}(H^A)_{\sigma_j}(x_{j+1},\tau_{j+1})
\Bigr)^2\Bigr)^{2^{k-j-1}}.$$

Before we continue, let us briefly see what principle was used when we
decided how to apply Cauchy-Schwarz. The idea was to take all terms that
did not depend on $x_{j+1}$ out to the left of $x_{j+1}$, except that
each time we took out a $(g^A)_{\sigma_j}$ or an $(H^A)_{\sigma_j}$, 
we left an $(H^A)_{\sigma_j}$ behind, exploiting the fact that 
$(g^A)_{\sigma_j}(H^A)_{\sigma_j}=(g^A)_{\sigma_j}$ and
$(H^A)_{\sigma_j}(H^A)_{\sigma_j}=(H^A)_{\sigma_j}$. In this way,
we extracted maximum information from the Cauchy-Schwarz inequality.

Since each $g^A$ is an $A$-function supported in $\ch$, and it 
maps to $[-1,1]$, and since each $H^A$ takes values $0$ or $1$,
we will not decrease the first term in the product
if we replace it by 
$$\Bigl(\E_{\tau_{j+1}}\prod_{[j]\subset A\atop {j+1\notin A\atop |A|<k}}
(H^A)_{\sigma_j}(x_{j+1},\tau_{j+1})
\prod_{[j]\not\subset A\atop {j+1\notin A\atop|A|<k}}(H^A)_{\sigma_j}
(x_{j+1},\tau_{j+1})\Bigr)^{2^{k-j-1}},$$
which we can write more succinctly as 
$$\Bigl(\E_{\tau_{j+1}}\prod_{j+1\notin A\atop |A|<k}
(H^A)_{\sigma_j}(\tau_j)\Bigr)^{2^{k-j-1}}.$$
To deal with the second term, we first have to expand out the
square, which in our notation is rather simple: we obtain
$$\Bigl(\E_{x_{j+1},x_{j+1}'}\E_{\tau_{j+1}}f_{\sigma_{j+1}}(\tau_{j+1})
\prod_{[j+1]\subset A}(g^A)_{\sigma_{j+1}}(\tau_{j+1})
\prod_{[j+1]\not\subset A\atop |A|<k}(H^A)_{\sigma_{j+1}}(\tau_{j+1})
\Bigr)^{2^{k-j-1}}.$$

We now apply H\"older's inequality. This time we take to the left of
the expectation over $\tau_{j+1}$ all terms that have no dependence on
$\tau_{j+1}$, again leaving behind the corresponding
$(H^A)_{\sigma_{j+1}}$ terms as we do so. The one exception is that,
for convenience only, we do not take the term $(g^A)_{\sigma_{j+1}}$
to the left when $A=[j+1]$, but instead take out
$(H^A)_{\sigma_{j+1}}$ in this case.  The result is that the last
quantity is bounded above by the product of
$$\Bigl(\E_{x_{j+1},x_{j+1}'}
\prod_{A\subset[j+1]\atop|A|<k}
H^A_{\sigma_{j+1}}\Bigr)^{2^{k-j-1}-1}$$
and
$$\E_{x_{j+1},x_{j+1}'}\Bigl(\E_{\tau_{j+1}}f_{\sigma_{j+1}}(\tau_{j+1})
\prod_{[j+1]\subset A}(g^A)_{\sigma_{j+1}}(\tau_{j+1})
\prod_{[j+1]\not\subset A\atop |A|<k}(H^A)_{\sigma_{j+1}}(\tau_{j+1})
\Bigr)^{2^{k-j-1}}.$$
 
These calculations have given us the expression we started with,
inside an expectation, with $j$ replaced by $j+1$. We must therefore 
check that we also have a chain $\ck_{j+1}$ with the right properties. 
Looking back at the various brackets we have discarded, this tells us 
that we want to rewrite the expression
$$\E_{\omega_j}\Bigl(\prod_{A\in\ck_j}H^A(\omega_j)\Bigr)
\Bigl(\E_{\tau_{j+1}}\prod_{j+1\notin A\atop |A|<k}
(H^A)_{\sigma_j}(\tau_j)\Bigr)^{2^{k-j-1}}
\Bigl(\E_{x_{j+1},x_{j+1}'}
\prod_{A\subset[j+1]\atop|A|<k}
H^A_{\sigma_{j+1}}\Bigr)^{2^{k-j-1}-1}$$
as 
$$\E_{\omega_{j+1}}\Bigl(\prod_{A\in\ck_{j+1}}H^A(\omega_{j+1})\Bigr)$$
for a chain $\ck_{j+1}$ with properties analogous to those of $\ck_j$.

There is a slight abuse of notation above, because after our applications
of the Cauchy-Schwarz and H\"older inequalities we have ended up
overusing $\tau_{j+1}$, $x_{j+1}$ and $x_{j+1}'$. But we can cure this
by renaming the variables in the expression we wish to
rewrite. Indeed, since we are raising the expectation over
$\tau_{j+1}=(x_{j+2},\dots,x_t)$ to the power $2^{k-j-1}$, let us 
introduce $2^{k-j-1}$ new variables for each variable included in 
$\tau_{j+1}$. More precisely, let us choose a set $U$ of cardinality
$2^{k-j-1}$ that is disjoint from $U_i^j$ for every $i$ between $j+1$
and $t$ and replace $V_i^j=\{i\}\times U_i^j$ by $\{i\}\times(U_i^j\cup U)$.
We can then expand out the second bracket as an expectation over the 
variables $x_1,x_1',\dots,x_j,x_j'$ and $x_i^u$ with $i\ge j+2$ and 
$u\in U$ of the product of all expressions of the form 
$(H^A)_{\sigma_j}(\tau_j^u)$, where $\tau_j^u=(x_{j+1}^u,\dots,x_t^u)$.
(In fact, there is no dependence on $x_{j+1}^u$, but we add the 
variables anyway so that it looks slightly nicer.) 

In a similar way, we can expand out the third bracket and introduce a
further $2(2^{k-j-1}-1)$ new variables into $V_{j+1}^j$.  When we do
these expansions, we end up writing the expression in the desired form
for some set-system $\ck_{j+1}$. It is not hard to see that
$\ck_{j+1}$ is a chain, so it remains to prove that it contains the
right number of sets of each index. 

Let $\g$ be the usual projection $(i,h)\mapsto i$. We need to prove
that each set $A\in\cj$ of cardinality less than $k$ has exactly $2^k$
preimages under $\g$ in $\ck_{j+1}$. We consider various cases.

First, if $A$ is a subset of $[j]$, then $\ck_j$ (which we
can think of as a chain defined on the vertex sets of $\ck_{j+1}$)
already contains $2^k-2^{|A|}$ preimages of $A$. Since the additional 
vertices $(i,u)$ do not project into $[j]$, we do not
create any new preimages in $\ck_{j+1}$.

Now suppose that $A$ is a subset of $[j+1]$ that contains
$j+1$. Then $A\not\subset[j]$ so the number of preimages of 
$A$ in $\ck_j$ is $2^k-2^{k-j+|A\cap[j]|}$. No new preimages come
from the second bracket, since that involves only sets that do not
include $j+1$, while from the third bracket we obtain 
$(2^{|A\cap[j+1]|})(2^{k-j-1}-1)$ preimages. But 
$2^{k-j-1+|A\cap[j+1]|}=2^{k-j+|A\cap[j]|}$ in this case, so
the total number of preimages is $2^k-2^{|A\cap[j+1]|}=2^k-2^{|A|}$.

Next, suppose that $A\not\subset[j+1]$ and $j+1\in A$. Then 
$\ck_j$ contains $2^k-2^{k-j+|A\cap[j]|}$ preimages of $A$ and the
second and third brackets do not contribute any. Since
$k-j+|A\cap[j]|=k-j-1+|A\cap[j+1]|$, the total number of preimages
is $2^k-2^{k-j-1+|A\cap[j+1]|}$, as we want.

Finally, suppose that $A\not\subset[j+1]$ and $j+1\notin A$. 
In that case, $\ck_j$ contains $2^k-2^{k-j+|A\cap[j]|}$ preimages,
the third bracket contributes none, and the second bracket contributes
$2^{|A\cap[j]|}2^{k-j-1}=2^{k-j-1+|A\cap[j]|}$ preimages. Thus, the
total number of preimages is $2^k-2^{k-j-1+|A\cap[j]|}$, which equals
$2^k-2^{k-j-1+|A\cap[j+1]|}$.

This completes the proof of the inductive step. All that remains is
the simple task of checking that the case $j=k$ of the induction 
is the statement that we wish to prove. But when $j=k$, we have the
upper bound
$$\E_{\omega_k}\Bigl(\prod_{A\in\ck_k}H^A(\omega_k)\Bigr)
\Bigl(\E_{\tau_k}f_{\sigma_k}(\tau_k)
\prod_{[k]\subset A}(g^A)_{\sigma_k}(\tau_k)
\prod_{[k]\not\subset A\atop |A|<k}(H^A)_{\sigma_k}(\tau_k)\Bigr)^{2^{k-k}}.$$
The most obvious simplification is $1$ for $2^{k-k}$. Since $\cj$ does
not contain the set $[k]$, the first product in the second bracket
disappears.  This gives us the upper bound
$$\E_{\omega_k,\tau_k}f_{\sigma_k}(\tau_k)\prod_{A\in\ck_k}H^A(\omega_k)
\prod_{|A|<k}(H^A)_{\sigma_k}(\tau_k).$$
Writing $\omega$ for $(\omega_k,\tau_k)$ and letting $\ck$ be the
union of the sets in $\ck_k$ and the sets implied by the second
product (we will say what these are in a moment), we can write
this as
$$\E_\omega f_\sigma(\omega)\prod_{A\in\ck}H^A(\omega)$$
as required.

We still need to check that $\ck$ contains precisely $2^k$ preimages
of each set $A\in\cj$ of cardinality less than $k$. Let us therefore
be slightly more explicit about the ``sets implied by the second
product.'' A function $(H^A)_{\sigma_k}(\tau_k)$ is a product of
functions of the form $H^A(x_1^{\e_1},\dots,x_k^{\e_k},\tau_k)$. But
$H^A$ depends only on the variables in $A$, so the number of distinct
functions in the product is $2^{|A\cap[k]|}$, and thus the number of
preimages of $A$ in $\ck$ that come from the second product is
$2^{|A\cap[k]|}$. But when $j=k$, the number of preimages in $\ck_k$
is $2^k-2^{|A\cap[k]|}$, whether or not $A$ is a subset of $[k]$.
Therefore, for each set $C\subset\{1,2,\dots,r\}$ of cardinality less
than $k$, the chain $\ck$ contains precisely $2^k$ sets of index $C$
for each set $A\in\cj$ of index $C$, as claimed.
\hfill $\square$ \bigskip

As we shall see in the next section, the fact that the sets in $\ck$
have cardinality at most $k-1$ allows us to use Lemma 4.1 inside another
induction (in fact, a double induction). This corresponds to the 
second part of \S 2, where we replaced functions such as $G(x,t)$
by constant functions $\d_{XT}$. This time the functions we shall 
replace are functions of the form $H^A$ with $A\in\ck$. 
\bigskip

\noindent {\bf \S 5. A counting lemma for quasirandom chains.}
\medskip

Just before we prove our main result, we isolate a simple statement
that is needed in the proof and that helps to explain some of our 
choices in the definition of $(\e,\cj,k)$-quasirandom chains. For
convenience, we briefly recall the definition here. We constructed
a sequence $\e_k,\e_{k-1},\dots,\e_1$ by letting $\e_k=\e$ and
$$\e_{k-j}=2^{-jk-1}|\cj|^{-1}\Bigl(\e_{k-j+1}
\prod_{A\in\cj\atop|A|\ge k-j+1}\d_A\Bigr)^{2^{jk}}$$
when $j\ge 1$. We also defined $\eta_{k-j}$ by the formula
$$\eta_{k-j}=(1/2)
\Bigl(\e_{k-j}\prod_{A\in\cj\atop|A|\ge k-j}\d_A\Bigr)^{2^{k(j+1)}}$$
for each $j$. Finally, we declared $\ch$ to be $(\e,\cj,k)$-quasirandom
if, for every $A\in\cj$ of size $j\le k$, the hypergraph $H(A)$
was $\eta_j$-quasirandom relative to $H_*(A)$.

These parameters are chosen in order to satisfy some assumptions 
required in the inductive step of Theorem 5.2 below. The next lemma 
establishes that they do indeed satisfy them.

\proclaim Lemma {5.1}. Let $\cj$ and $\ch$ be chains and suppose
that $\ch$ is $(\e,\cj,k)$-quasirandom. Let $\ck$ be a chain with
the same vertex set as that of $\cj$, and suppose that there is
a homomorphism from $\ck$ to $\cj$ such that each set in $\cj$ 
has at most $2^k$ preimages. Let $\e_k,\e_{k-1},\dots,\e_1$ be the 
sequence defined above. Then $\ch$ is $(\e_{k-1},\ck,k-1)$-quasirandom.

\Proof  Let $\theta=\e_{k-1}$ and define a sequence 
$\theta_{k-1},\theta_{k-2},\dots$ by taking $\theta_{k-1}=\theta$ and
$$\theta_{k-1-j}=2^{-j(k-1)-1}|\ck|^{-1}\Bigl(\theta_{k-j}
\prod_{A\in\ck\atop|A|\ge k-j}\d_A\Bigr)^{2^{j(k-1)}}.$$
Suppose that $\theta_{k-j}\ge\e_{k-j}$. We also know that
$|\ck|^{-1}\ge 2^{-k}|\cj|^{-1}$ and that 
$$\prod_{A\in\ck\atop|A|\ge k-j}\d_A\ge
\Bigl(\prod_{A\in\cj\atop|A|\ge k-j}\d_A\Bigr)^{2^k}.$$
It follows that
$$\theta_{k-1-j}\ge 2^{-jk-1}|\cj|^{-1}
\Bigl(\e_{k-j}\prod_{A\in\cj\atop|A|\ge k-j}\d_A\Bigr)^{2^{jk}}=\e_{k-j}.$$
Therefore by induction $\theta_j\ge\e_j$ for every $j$. 

Now let $j$ be an integer between $0$ and $k-1$. Then
$$\eqalign{\eta_{k-1-j}=\eta_{k-(j+1)}&
=(1/2)\Bigl(\e_{k-(j+1)}
\prod_{A\in\cj\atop|A|\ge k-(j+1)}\d_A\Bigr)^{2^{k(j+2)}}\cr
&\le(1/2)\Bigl(\theta_{k-(j+1)}
\Bigl(\prod_{A\in\cj\atop|A|\ge k-(j+1)}\d_A\Bigr)^{2^k}\Bigr)^{2^{k(j+1)}}\cr
&\le(1/2)\Bigl(\theta_{k-(j+1)}
\prod_{A\in\ck\atop|A|\ge k-(j+1)}\d_A\Bigr)^{2^{k(j+1)}}\cr
&\le(1/2)\Bigl(\theta_{k-1-j}
\prod_{A\in\ck\atop|A|\ge k-1-j}\d_A\Bigr)^{2^{(k-1)(j+1)}}\ .\cr}$$
This is the formula for $\eta_{k-j}$ except that $k$ has been replaced
by $k-1$, $\cj$ by $\ck$, and $\e_{k-j}$ by $\theta_{k-1-j}$. It 
follows that $\ch$ is $(\e_{k-1},\ck,k-1)$-quasirandom, as claimed.
\hfill $\square$ \bigskip

In the next theorem and its proof, we shall discuss two chains
$\cj$ and $\ch$, and borrow notation from the previous section
without redefining it. For example, $\tau$ is once again a
sequence $(\seq x t)$ that enumerates variables that are indexed
by the vertices of $\cj$. Eventually, we will be interested in the
case where every function $g^A$ is just $H^A$, but this more
general statement is needed for an inductive argument to work,
and is also of some interest in its own right.

\proclaim Theorem {5.2}. Let $\cj$ and $\ch$ be $r$-partite chains 
as described at the beginning of the previous section. Let $\cj_1$ be
a subchain of $\cj$ and for each $A\in\cj_1$ let $g^A$ be an
$A$-function supported in $\ch$. Suppose that the maximum cardinality
of any set in $\cj\setminus\cj_1$ is $k$ and that $\ch$ is 
$(\e,\cj,k)$-quasirandom. Then
$$\Bigl|\E_\tau\prod_{A\in\cj_1}g^A(\tau)\prod_{A\in\cj\setminus\cj_1}
H^A(\tau)-
\E_\tau\prod_{A\in\cj_1}g^A(\tau)\prod_{A\in\cj\setminus\cj_1}\d_A\Bigr|
\le\e|\cj\setminus\cj_1|\prod_{A\in\cj}\d_A\ .$$

\Proof  This result tells us that we can replace the functions
$H^A$ in the quantity
$\E_\tau\prod_{A\in\cj_1}g^A(\tau)\prod_{A\in\cj\setminus\cj_1}
H^A(\tau)$ by their relative densities $\d_A$ without changing the
quantity by too much. This is proved by two levels of induction, for
the following reason. First of all, we do our replacements one by
one, and this leads to an induction on the cardinality
of $\cj\setminus\cj_1$. However, in order to establish an upper bound for
the error introduced when we make a replacement, we use our main lemma,
Lemma 4.1, which results in a similar expression to the one we were
initially trying to bound, but with new chains $\ck$ and $\ck_1$. These 
chains are considerably bigger than $\cj$ and $\cj_1$, but the largest set in
$\ck\setminus\ck_1$ is smaller than the largest set in
$\cj\setminus\cj_1$, so we can use induction on $k$ to replace the
error term itself by a quantity that will turn out to be small as a
direct consequence of the quasirandomness of the chain~$\ch$.

Let us therefore choose a maximal set $A_0$ in $\cj\setminus\cj_1$
and try to replace $H^{A_0}(\tau)$ by $\d_{A_0}$ in the expression 
$\E_\tau\prod_{A\in\cj_1}g^A(\tau)\prod_{A\in\cj\setminus\cj_1}
H^A(\tau)$ while introducing only a small error. Letting 
$\cj_0=\cj\setminus\{A_0\}$, the difference between the original
expression and the new expression is
$$\E_\tau f(\tau)\prod_{A\in\cj_1}g^A(\tau)
\prod_{A\in\cj_0\setminus\cj_1}H^A(\tau),$$
where $f$ is the $A_0$-function defined by 
$f(\tau)=(H^{A_0}(\tau)-\d_{A_0})\prod_{A\subsetneq A_0}H^A(\tau)$.
(This function was first defined near the end of subsection 3.6: 
in the notation of this section it equals $1-\d_{A_0}$ if 
$\tau(A_0)\in H(A_0)$, $-\d_{A_0}$ if 
$\tau(A_0)\in H_*(A_0)\setminus H_{A_0}$, and zero otherwise.)

Without loss of generality, we may assume that $A_0$ is the
set $\{1,2,\dots,k\}$. Let us therefore apply Lemma 4.1 to
this function $f$ and to the chain $\cj_0$. It yields for us
an $r$-partite $(k-1)$-chain $\ck'$ and a homomorphism $\g$ 
from $\ck'$ to $\cj_0$ such that every set in $\cj_0$ of
cardinality less than $k$ has $2^k$ preimages, and such that
we have the inequality
$$\Bigl(\E_\tau f(\tau)\prod_{A\in\cj_1}g^A(\tau)
\prod_{A\in\cj_0\setminus\cj_1}H^A(\tau)\Bigr)^{2^k}
\le \E_\omega f_{\sigma}(\omega)\prod_{B\in\ck'}H^B(\omega).$$

Recall that $f(\sigma)$ is the product of $f(\omega(A))$ over all sets
$A$ of the form $\{(1,\e_1),\dots,(k,\e_k)\}$.  Let $\ck_1$ be the
chain of all subsets of such sets and let $\ck=\ck_1\cup\ck'$. Then
the largest set in $\ck\setminus\ck_1$ has size at most
$k-1$. Moreover, by Lemma 5.1, $\ch$ is $(\e_{k-1},\ck,k-1)$-quasirandom.
Therefore, by induction on $k$, we know that the right-hand
side of the above inequality differs from $\E_\sigma
f_{\sigma}\prod_{A\in\ck\setminus\ck_1}\d_{A}$ by at most
$\e_{k-1}|\ck\setminus\ck_1|\prod_{A\in\ck}\d_A$.

This is at most $\e_{k-1}|\ck\setminus\ck_1|\prod_{A\in\ck'}\d_A$,
which is equal to $\e_{k-1}|\ck\setminus\ck_1|
\Bigl(\prod_{A\in\cj_0\atop|A|<k}\d_A\Bigr)^{2^k}$. But
$|\ck\setminus\ck_1|\le|\ck|\le 2^k|\cj|$ and
$2^{k+1}\e_{k-1}|\cj|\le
\Bigl(\e_k\prod_{A\in\cj\atop|A|\ge k}\d_A\Bigr)^{2^k}$,
so this is at most $(1/2)\Bigl(\e_k\prod_{A\in\cj}\d_A\Bigr)^{2^k}$.

As for $\E_\sigma 
f_{\sigma}\prod_{A\in\ck\setminus\ck_1}\d_{A}$, it is equal (by
definition) to $\oct(f)\prod_{A\in\ck\setminus\ck_1}\d_{A}$. By
hypothesis, $f$ is $\eta_k$-quasirandom, which means that
$\oct(f)\le\eta_k\prod_{A\in\ck_1}\d_A$. Since 
$\eta_k\le(1/2)\Bigl(\e_k\prod_{A\in\cj\atop|A|\ge k}\d_A\Bigr)^{2^k}$,
it follows that
$$\E_\sigma f_{\sigma}\prod_{A\in\ck\setminus\ck_1}\d_{A}\le
\eta_k\prod_{A\in\ck}\d_A\le\eta_k\prod_{A\in\ck'}\d_A
\le(1/2)\Bigl(\e_k\prod_{A\in\cj}\d_A\Bigr)^{2^k}.$$

Putting these two estimates together, we find that
$$\Bigl|\E_\tau f(\tau)\prod_{A\in\cj_1}g^A(\tau)
\prod_{A\in\cj_0\setminus\cj_1}H^A(\tau)\Bigr|\le\e_k\prod_{A\in\cj}\d_A.$$

Thus, returning to the beginning of the proof, we have shown that replacing 
$H^A$ by $\d_A$ for any maximal element of $\cj\setminus\cj_1$ results in 
an error of at most $\e_k\prod_{A\in\cj}\d_A$. Therefore the result follows
by induction on $|\cj\setminus\cj_1|$ and the triangle inequality (and
the fact that $\e_k=\e$).
\hfill $\square$ \bigskip

If we now consider the case when $\cj_1$ is empty, then we obtain
the following corollary, which is the counting lemma that we have
been aiming for.

\proclaim Corollary {5.3}. Let $\cj$ and $\ch$ be $r$-partite chains with
vertex sets $E_1\cup\dots\cup E_r$ and $X_1\cup\dots\cup X_r$,
respectively. Let $k$ be the size of the largest set in $\cj$
and suppose that $\ch$ is $(\e/|\cj|,\cj,k)$-quasirandom. Let $\tau$
be a random $r$-partite map from $E_1\cup\dots\cup E_r$ to 
$X_1\cup\dots\cup X_r$. Then
$$\Bigl|\P[\tau\in\Hom(\cj,\ch)]-\prod_{A\in\cj}\d_A\Bigr|
\le\e\prod_{A\in\cj}\d_A\ .\eqno{\square}$$

In less precise terms, this says that if $\cj$ is a small $r$-partite
chain and $\ch$ is a sufficiently quasirandom $r$-partite chain, then
a random $r$-partite map from the vertices of $\cj$ to the vertices of
$\ch$ will be a homomorphism with approximately the probability that
you would expect if $\ch$ was a random chain with the given relative
densities.
\bigskip

\noindent{\bf 6. Local increases in mean-square density.}
\medskip

All known proofs of Szemer\'edi's theorem use (explicitly or implicitly) an 
approach of the following kind. Given a dense set that fails to be quasirandom
in some appropriate sense, one can identify imbalances in the set that allow 
one to divide it into pieces that ``improve'' in some way, on average at least,
on the set itself. One then iterates this argument until one reaches sets
that {\it are} quasirandom. At that point one uses some kind of counting lemma
to prove that they contain an arithmetic progression of length $k$.

This proof is no exception. We have defined a notion of quasirandomness
and proved a counting lemma for it. Now we must see what happens when
some parts of a chain are {\it not} relatively quasirandom. We shall 
end up proving a {\it regularity lemma}, which says, roughly speaking,
that any dense chain can be divided up into a bounded number of pieces,
almost all of which are quasirandom. This generalizes Szemer\'edi's
regularity lemma for graphs (which formed part of his proof of his
theorem on arithmetic progressions).

Given a dense graph $G$ and a positive real number $\e$, Szemer\'edi's
regularity lemma asserts that the vertices of $G$ can be partitioned
into $K$ classes of roughly equal size, with $K$ bounded above by a
function of $\e$ only, in such a way that, proportionately speaking,
at least $1-\e$ of the bipartite graphs spanned by two of these
classes are $\e$-regular. (One can insist that $K$ is much bigger
than $\e^{-1}$, so it is not necessary to worry about the case where the
two classes are equal. Or it can be neater to say that two equal
classes form a ``regular pair'' if they span a quasirandom graph.)

Very roughly, the proof is as follows. Suppose you have a graph $G$
and a partition of its vertex set. Then either this partition will do
or there are many pairs of cells from the partition that give rise to
induced bipartite subgraphs of $G$ that are not $\e$-quasirandom. If
$X$ and $Y$ are two disjoint sets of vertices, write $G(X,Y)$ for the
corresponding induced bipartite subgraph of $G$. Suppose that $X$ and
$Y$ are two cells of the partition, for which $G(X,Y)$ is not
$\e$-regular. Then there are large subsets $X(0)\subset X$ and
$Y(0)\subset Y$ for which the density of $G(X(0),Y(0))$ is
substantially different from that of $G(X,Y)$. Letting
$X(1)=X\setminus X(0)$ and $Y(1)=Y\setminus Y(0)$, we have obtained
partitions of $X$ and $Y$ into two sets each, in such a way that the
densities of the graphs $G(X(i),Y(j))$ are not almost all
approximately the same as that of $G(X,Y)$. One can then define an
appropriately weighted average of the squares of these four densities
and show that this average is greater than the square of the density
of $G(X,Y)$. Let us call this {\it stage one} of the argument, the
stage where we identify a ``local'' increase in mean-square density.

It remains to turn these local increases into a global increase.
This, which we shall call {\it stage two}, is quite simple.
Denote the cells of the original partition by $\seq X k$. For
each pair $(X_i,X_j)$ that fails to be $\e$-regular, use the
above argument to partition $X_i$ into two sets $X_{ij}(0)$ and 
$X_{ij}(1)$, and to partition $X_j$ into two sets $X_{ji}(0)$ 
and $X_{ji}(1)$. Then for each $i$ find a partition of $X_i$
that refines all the partitions $\{X_{ij}(0),X_{ij}(1)\}$.
The result is a partition into $m\le k.2^k$ sets $\seq Y m$ that 
refines the partition $\{\seq X k\}$. It can be shown that the 
average of the squares of the densities $G(Y_i,Y_j)$, again,
with appropriate weights, is significantly greater than it
was for the partition $\{\seq X k\}$. Therefore, if one iterates
the procedure, the iteration must terminate after a number of
steps that can be bounded in terms of $\e$. It can terminate
only if almost all the graphs $G(X_i,X_j)$ are quasirandom, so
the result is proved.

We have given this sketch since our generalized regularity lemma
will be proved in a similar way. There are two main differences.
First, it is an unfortunate fact of life that, when one is dealing
with $k$-chains rather than graphs, simple arguments have to
be expressed in terminology that can obscure their simplicity.
For example, even defining the appropriate notion of a ``partition'' 
of a chain is somewhat complicated. Thus, stage two of our argument,
although it is an ``obvious'' generalization of stage two of the proof
of the usual regularity lemma, is noticeably more complicated to
write down.

A more fundamental difference, however, is that our stage one is not
completely straightforward, and here the difference is mathematical
rather than merely notational. The reason is that we do not generalize
Szemer\'edi's regularity lemma as it is stated above, but rather a
simple variant of it where rather than obtaining $\e$-regular pairs we
obtain $\e$-{\it quasirandom} pairs. For dense bipartite graphs, these
two notions are equivalent (give or take changes in $\e$), but when
one generalizes them to hypergraphs that live in sparse chains they
diverge in a significant way. Some hint of this can already be seen
above. It is true by definition that if a pair $G(X,Y)$ is not
$\e$-regular, then there are large subsets $X(0)\subset X$ and
$Y(0)\subset Y$ for which the density of $G(X(0),Y(0))$ is
substantially different from that of $G(X,Y)$. However, if we assume
instead that $G(X,Y)$ is not $\e$-quasirandom, then there is something
to prove. The proof is very simple in the dense case, and even in the
sparse case, but in the latter it yields sets $X(0)$ and $Y(0)$ that
are very small. As a result, we have to work significantly harder in
order to obtain a partition with a good enough local increase in mean-square
density. Roughly speaking, our approach will be to find many pairs of
such sets, and build a partition out of those. For this to work it is
important that the pairs are sufficiently spread out: the detailed argument 
will occupy the rest of the section.

Incidentally, the last paragraph describes the main difference between
our approach and that of Nagle, R\"odl, Schacht and Skokan. Their
definitions generalize that of $\e$-regularity of bipartite graphs,
so stage 1 of the proof of the regularity lemma is easier for them.
However, they have to pay for this when they prove their counting 
lemma: $\e$-regularity is a weaker property than $\e$-quasirandomness,
so if you use it as your basic definition then it is easier to deduce
facts about objects that are {\it not} $\e$-regular but harder to
deduce facts about objects that {\it are} $\e$-regular.

We shall now work towards our stage one, which will be  
Lemma 6.3 below. To begin with, let us say what we mean by the
mean-square density of a function with respect to a partition. Let $U$
be a set of size $n$, let $f:U\ra\R$ and let $\seq B r$ be sets that
form a partition of $U$. Then the {\it mean-square density} of $f$
{\it with respect to} the partition $\{\seq B r\}$ is
$$\sm i r{|B_i|\over n}\Bigl(\E_{x\in B_i}f(x)\Bigr)^2\ .$$
If we write $\b_i$ for $|B_i|/n$ (which it is helpful to think of
as the probability that a random $x\in U$ is an element of $B_i$)
and $\d_i$ for $\E_{x\in B_i}f(x)$ (that is, the expectation,
or ``density'', of $f$ in $B_i$) then this sum is $\sm i r\b_i\d_i^2$,
the weighted average of the squared densities $\d_i^2$, with respect to 
the obvious system of weights $\b_i$.

The following two simple lemmas are very slight modifications of
lemmas in [G2]. The first is our main tool, while the second is more
of a technical trick that will be used in Lemma 6.3.

\proclaim Lemma {6.1}. Let $U$ be a finite set and let $f$ and $g$
be functions from $U$ to the interval $[-1,1]$. Let $\seq B r$ be a 
partition of $U$ and suppose that $g$ is constant on each $B_i$.
Then the mean-square density of $f$ with respect to the partition
$\seq B r$ is at least $\sp{f,g}^2/\|g\|_2^2$. 

\Proof  For each $j$ let $a_j$ be the value taken by $g$ on the set $B_j$.
Then, by the Cauchy-Schwarz inequality,
$$\eqalign{\sp{f,g}^2&=\Bigl(\sum_ja_j\b_j\E_{x\in B_j}f(x)\Bigr)^2\cr
&\le\Bigl(\sum_j\b_ja_j^2\Bigr)
\Bigl(\sum_j\b_j\Bigl(\E_{x\in B_j}f(x)\Bigr)^2\Bigr)\ .\cr}$$
The first part of the product is $\|g\|_2^2$ and the second is
the mean-square density of $f$, from which the lemma follows. 
\hfill $\square$ \bigskip

In the next lemma, $\E_iv_i$ and $\E_iw_i$ mean the obvious thing: 
they are $n^{-1}\sm i nv_i$ and $n^{-1}\sm i nw_i$, respectively. 

\proclaim Lemma {6.2}. Let $n$ be a positive integer, let 
$0<\d<1$ and let $r$ be an integer greater than or equal to
$\d^{-1}$. Let $\seq v n$ be vectors in a Hilbert space such that
$\|v_i\|^2\le 1$ for each $i$ and such that $\left\|\E_i
v_i\right\|^2\le\d$.  Let $r$ vectors $\seq w r$ be chosen
uniformly and independently from the $v_i$.  (To be precise, for each
$w_j$ an index $i$ is chosen randomly between 1 and $n$ and $w_j$ is
set equal to $v_i$.) Then the expectation of $\left\|\E_j
w_j\right\|^2$ is at most~$2\d$.

\Proof  The expectation of $\left\|\E_j w_j\right\|^2$ is the expectation
of $\E_{i,j}\sp{w_i,w_j}$. If $i\ne j$ then the expectation of $\sp{w_i,w_j}$
is $\left\|\E_iv_i\right\|^2$ which, by hypothesis, is at most $\d$.
If $i=j$, then $\sp{w_i,w_j}$ is at most $1$, again by hypothesis. Therefore,
the expectation we are trying to bound is at most $r^{-2}(\d r(r-1)+r)$. Since
$\d r\ge 1$, this is at most $2\d$, as claimed. \hfill $\square$ \bigskip

Before we state the main result of this section, we need two
definitions.  The first is of a chain $\cd$ that we shall call a {\it
double octahedron}. We use this name for conciseness even though it is
slightly misleading: in fact, $\cd$ is the $(k-1)$-skeleton of a chain
formed from two $k$-dimensional octahedra by identifying a face from
one with the corresponding face from the other. To put this more
formally, take the vertex set of $\cd$ to be the set
$[k]\times\{0,1,2\}$. For each $i$ between 1 and $k$ let $V_i$ be the
set $\{i\}\times\{0,1,2\}$ and for $j=0,1,2$ let $B_j$ be the set
$[k]\times\{j\}$. The edges of $\cd$ are all sets $B$ of cardinality
at most $k-1$ such that $|B\cap V_i|\le 1$ for every $i$ and at least
one of $B\cap B_1$ and $B\cap B_2$ is empty. (The two octahedra in
question are $O_1$ and $O_2$, where $O_j$ consists of all sets
$B\subset B_0\cup B_j$ such that $|B\cap V_i|\le 1$ for every $i$.)

Notice that if $A\subset[k]$ is a set of size at most $k-1$ then
the number of edges in $\cd$ of index $A$ is $2^{|A|+1}-1$,
since there are $2^{|A|}$ edges from each octahedron and
one, namely $A\times\{0\}$, which is common to both. 

For the second definition, suppose we have a $k$-partite $(k-1)$-chain
$\ch$ with vertex sets $\seq X k$. Recall from \S 2 that $H_*([k])$
is the collection of all sets $A$ such that $|A\cap X_i|=1$ for every
$i$ and such that every proper subset of $A$ belongs to $\ch$. For
this second condition to hold it is enough for $C$ to be an edge of
$\ch$ whenever $C\subset A$ and $|C|=k-1$. Let $H$ be the $k$-partite
$(k-1)$-uniform hypergraph consisting of all edges of $\ch$ of size
$k-1$. For $1\le i\le k$ let $H_i$ be the $(k-1)$-partite
subhypergraph of $H$ consisting of all edges of $H$ that have empty
intersection with $X_i$. We shall call the hypergraphs $H_i$ the {\it
parts} of $H$. Each set $A\in H_*([k])$ has $k$ subsets of size
$k-1$. Each part $H_i$ of $H$ contains exactly one of these subsets,
namely $A\setminus X_i$.

Suppose that each $H_i$ is partitioned into subhypergraphs
$H_{i1},\dots,H_{ir_i}$. These partitions give rise
to an equivalence relation $\sim$ on $H_*([k])$: we say
that $A\sim A'$ if, for each $i\le k$, the sets 
$A\setminus X_i$ and $A'\setminus X_i$ belong to the same 
cell $H_{ij}$ of the partition of $H_i$. The corresponding
partition will be called the {\it induced partition} of
$H_*([k])$. 

\proclaim Lemma {6.3}. Let $\ch$ be a $k$-partite $(k-1)$-chain with vertex
sets $\seq X k$, let $\cd$ be the double octahedron, let
$\d=\prod_{A\in\cd}\d_A$ and let $r\ge\d^{-1}$ be a positive integer.
Suppose that $\e\le|\cd|^{-1}$, that $\ch$ is
$(\e,\cd,k-1)$-quasirandom and that $f:H_*([k])\ra[-1,1]$ is a
function that is not $\eta$-quasirandom relative to $\ch$. Let $H$ be
the set of all edges of $\ch$ of size $k-1$ and let $\seq H k$ be the
$k$ parts of $H$. Then there are partitions of the $H_i$ into at most
$3^r$ sets each such that the mean-square density of $f$ with respect to
the induced partition of $H_*([k])$ is at least $\eta^2/32$.

We shall prove Lemma 6.3 in stages, by means of some intermediate
lemmas (Lemmas 6.4-6.7 below). Since these lemmas form part of a
larger proof, we shall not state each one in full: rather, if we have
already introduced notation such as names for various functions we
shall feel to use it again without redefining it.

But before we get on to the subsidiary lemmas, let us examine our main
hypothesis, that $f$ is not $\eta$-quasirandom relative to $\ch$. For
each $i\le k$ let $U_i=\{i\}\times\{0,1\}$ (so $U_i$ consists of the
``first two'' of the three elements of $V_i$). As in \S 3, let $\cb$
be the $k$-partite $k$-uniform hypergraph consisting of all sets
$B\subset U_1\cup\dots\cup U_k$ such that $|B\cap U_i|=1$ for every
$i$, let $\ck$ be the chain of all sets $C$ that are proper subsets of
some $B\in\cb$ and let $\Omega$ be the set of all $k$-partite maps
from $U_1\cup\dots\cup U_k$ to $X_1\cup\dots\cup X_k$.  Then to say
that $f$ is not $\eta$-quasirandom relative to $\ch$ is to say that
$$\oct(f)=\E_{\omega\in\Omega}\prod_{B\in\cb}f^B(\omega)
>\eta\prod_{A\in\ck}\d_A\ ,$$
where by $f^B(\omega)$ we mean $f(\omega(B))$ if 
$\omega(B)\in H_*([k])$ and 0 otherwise. 

Let $B_0$ and $B_1$ be as defined earlier, so that 
$U_1\cup\dots\cup U_k=B_0\cup B_1$. Let $\Phi$ and $\Psi$ be the
set of all $k$-partite maps from $B_0$ and $B_1$, respectively, 
to $X_1\cup\dots\cup X_k$. There is an obvious one-to-one
correspondence between $\Omega$ and $\Phi\times\Psi$: given 
any $\omega\in\Omega$, associate with it the pair $(\phi,\psi)$ where
$\phi$ and $\psi$ are the restrictions of $\omega$ to $B_0$ and
$B_1$. This procedure is invertible: given a pair $(\phi,\psi)$, 
define a $k$-partite map $\omega$ by setting $\omega(x)=\phi(x)$ 
if $x\in B_0$ and $\omega(x)=\psi(x)$ if $x\in B_1$. From now
on we shall identify $\Omega$ with $\Phi\times\Psi$ and freely
pass from one to the other. 

Let us split the product $\prod_{B\in\cb}f^B(\omega)$ into
two parts. We shall write $F(\omega)$ for $f^{B_0}(\omega)$
and $G(\omega)$ for $\prod_{B\in\cb,B\ne B_0}f^B(\omega)$. Now
if $\omega=(\phi,\psi)$ then $F(\omega)$ does not depend on $\psi$
(since it depends only on $\omega(B_0)=\phi(B_0)$). To emphasize
this, we shall write $G(\phi,\psi)$ for $G(\omega)$ and
$F(\phi)$ for $F(\omega)$. Our hypothesis now becomes
$$\E_{\phi\in\Phi}\E_{\psi\in\Psi}F(\phi)G(\phi,\psi)
>\eta\prod_{A\in\ck}\d_A\ . \eqno(*) $$

Let us see why this is useful. First, note that there is another 
obvious one-to-one correspondence, this time between $\Phi$ and 
$X_1\times\dots\times X_k$. It associates with a map $\phi\in\Phi$
the $k$-tuple $(\phi(1,0),\dots,\phi(k,0))$, and the inverse 
associates with a $k$-tuple $(\seq x k)\in\prod_{i=1}^kX_i$ 
the map $\phi:B_0\ra X_1\cup\dots\cup X_k$ that
takes $(i,0)$ to $x_i$ for each $i\le k$. Therefore, the 
function $F$ is basically another way of thinking about $f$.
The inequality above can be regarded as saying that, for
an average $\psi\in\Psi$, $F$ has a certain correlation with
the function $G_\psi:\phi\mapsto G(\phi,\psi)$. This is
significant, because the functions $G_\psi$ have a special
form, as the next lemma shows.

\proclaim Lemma {6.4}. Each function $G_\psi:\Phi\ra[-1,1]$ defined
above can be written as a product of $A$-functions over sets $A\subset
B_0$ of size $k-1$.

\Proof  By definition, 
$G_\psi(\phi)=\prod_{B\in\cb,B\ne B_0}f^B(\phi,\psi)$. Now
$f^B(\phi,\psi)$ depends on 
$(\phi,\psi)(B)=\phi(B\cap B_0)\cup\psi(B\cap B_1)$ only.
Therefore, if $\psi$ is fixed, $f^B(\phi,\psi)$ depends on 
$\phi(B\cap B_0)$ only. Thus, the function $\phi\mapsto f^B(\phi,\psi)$ 
is a $(B\cap B_0)$-function defined on $\Phi$. Since $B\ne B_0$,
$|B\cap B_0|\le k-1$. This proves that $G_\psi$ is a product of
$A$-functions over sets $A$ of size at most $k-1$. However, 
if $B\subset A$, then the product of a $B$-function with an
$A$-function is still an $A$-function. From this simple
observation it now follows that $G_\psi$ is a product of $A$-functions
over sets $A$ of size equal to $k-1$. \hfill $\square$ \bigskip

Our next task is to construct some new functions $E_\psi$ out of
the $G_\psi$ that have very similar properties but take values
$0$, $1$ and $-1$ only.

\proclaim Lemma {6.5}. If the inequality ($*$) holds, then there exist 
functions $E_\psi:\Phi\ra\{-1,0,1\}$, one for each $\psi\in\Psi$, with
the following properties. First, $E_\psi(\phi)$ is non-zero only if
$(\phi,\psi)\in\Hom(\ck,\ch)$. Second, each $E_\psi$ can be written as
a product of $\{-1,0,1\}$-valued $A$-functions over subsets $A\subset
B_0$ of size $k-1$. Third,
$$\E_{\phi\in\Phi}\E_{\psi\in\Psi}F(\phi)E_\psi(\phi)
>\eta\prod_{A\in\ck}\d_A\ .$$

\Proof  Let us fix $\psi\in\Psi$ and consider the function $G=G_\psi$.
By Lemma 6.4 we can write it as a product of $A$-functions, where
each $A$ in the product is a subset of $B_0$ of size $k-1$. There
are $k$ such sets, namely $\seq A k$, where for each $i$ we set
$A_i=B_0\setminus\{(i,0)\}$. So we can write 
$G(\phi)=\prod_{i=1}^kg_i(\phi)$ with $g_i$ an 
$A_i$-function for each $i$.

Now define an $A_i$-function $u_i:\Phi\ra\{-1,0,1\}$ randomly in 
the following natural way. Say that two maps $\phi$ and $\phi'$
are {\it equivalent} if $\phi(A_i)=\phi'(A_i)$ and choose one
map from each equivalence class. Let $\phi$ be one of
these representatives. If $g_i(\phi)\ge 0$ then let $u_i(\phi)$
equal 1 with probability $g_i(\phi)$ and $0$ with probability
$1-g_i(\phi)$. If $g_i(\phi)<0$ then let $u_i(\phi)$ equal -1
with probability $-g_i(\phi)$ and $0$ with probability 
$1+g_i(\phi)$. Then the expectation of $u_i(\phi)$ is 
$g_i(\phi)$. If $\phi'$ is equivalent to $\phi$ then let
$u_i(\phi')=u_i(\phi)$. 

Do the same for each equivalence class and make all the random
choices independently. Finally, for each $\phi\in\Phi$ let 
$E_\psi(\phi)=\prod_{i=1}^ku_i(\phi)$.

Now $E_\psi(\phi)$ can be non-zero only if $u_i(\phi)\ne 0$
for every $i$, and this is the case (with probability 1) only 
if $g_i(\phi)\ne 0$ for every $i$, and hence only if $G(\phi)\ne 0$.
We defined $G(\phi)$ to be 
$G_\psi(\phi)=\prod_{B\in\cb,B\ne B_0}f^B(\phi,\psi)$. 
But $f^B(\phi,\psi)=0$ unless $(\phi,\psi)(B)\in H_*([k])$, and this
is true only if $(\phi,\psi)(C)\in\ch$ for every proper subset
$C$ of $B$. Therefore this product is non-zero only if $(\phi,\psi)$ is 
a homomorphism from $\ck$ to $\ch$.

Since the choices of the different functions $u_i$ were made
independently and the expectation of $u_i(\phi)$ is $g_i(\phi)$, 
the expectation of $u_1(\phi)\dots u_k(\phi)$ is 
$g_1(\phi)\dots g_k(\phi)=G_\psi(\phi)$. Therefore, by linearity
of expectation, the expectation of $\E_\phi\E_\psi F(\phi)E_\psi(\phi)$
is $\E_\phi\E_\psi F(\phi)G_\psi(\phi)$, which we have assumed
to be at least $\eta\prod_{A\in\ck}\d_A$. It follows
that we can choose functions $E_\psi$ with the desired properties.
\hfill $\square$ \bigskip

\proclaim Lemma {6.6}. For each $\psi\in\Psi$ let $E_\psi$ be the
function constructed in Lemma 6.5, and let $\cd$ be the double
octahedron chain introduced before the statement of Lemma 6.3. Then 
$$\E_{\phi\in\Phi}\Bigl(\E_{\psi\in\Psi}E_\psi(\phi)\Bigr)^2
\le 2\prod_{A\in\cd}\d_A\ .$$

\Proof  The left-hand side of the inequality we wish to prove can 
be rewritten
$$\E_{\phi\in\Phi}\E_{\psi_1,\psi_2\in\Psi}
E_{\psi_1}(\phi)E_{\psi_2}(\phi)\ .$$
By Lemma 6.5, $E_{\psi_1}(\phi)E_{\psi_2}(\phi)$ is non-zero
if and only if $(\phi,\psi_1)$ and $(\phi,\psi_2)$ belong to
$\Hom(\ck,\ch)$. Therefore, this sum is at most the probability,
for a random triple $(\phi,\psi_1,\psi_2)\in\Phi\times\Psi^2$,
that both $(\phi,\psi_1)$ and $(\phi,\psi_2)$ belong to $\Hom(\ck,\ch)$. 

In order to estimate this probability, we shall apply the counting
lemma to the chain $\cd$. Every edge of $\cd$ is a proper subset of
either $B_0\cup B_1$ or $B_0\cup B_2$.  Let $\ck_1$ be the set of all
edges of the first kind and let $\ck_2$ be the set of all edges of the
second kind. Both $\ck_1$ and $\ck_2$ are chains and they intersect in
a chain that consists of all proper subsets of $B_0$. Moreover,
$\ck_1$ is essentially the same chain as $\ck$ (formally, it has
different vertex sets but the edges are the same). As for $\ck_2$, it
is isomorphic to $\ck$ in the following sense. Let $\g$ be the
bijection from $B_0\cup B_2$ to $B_0\cup B_1$ that takes $(i,0)$ to
$(i,0)$ and $(i,2)$ to $(i,1)$. Then $A$ is an edge of $\ck_2$ if and
only if $\g(A)$ is an edge of $\ck$.

Let $\Theta$ be the set of all $k$-partite functions
from $V_1\cup\dots\cup V_k$ (the vertex set of $\cd$) to
$X_1\cup\dots\cup X_k$. There is a one-to-one correspondence 
between $\Theta$ and $\Phi\times\Psi\times\Psi$ taking
$\theta\in\Theta$ to $(\phi,\psi_1,\psi_2\circ\gamma)$, where $\phi$,
$\psi_1$ and $\psi_2$ are the restrictions of $\theta$ to $B_0$, 
$B_1$ and $B_2$, respectively. Since $\cd=\ck_1\cup\ck_2$,
a map $\theta\in\Theta$ belongs to $\Hom(\cd,\ch)$ if and only 
if $(\phi,\psi_1)$ belongs to $\Hom(\ck_1,\ch)$ and 
$(\phi,\psi_2)$ belongs to $\Hom(\ck_2,\ch)$.
But this is true if and only if $(\phi,\psi_1)$ and
$(\phi,\psi_2\circ\gamma)$ belong to $\Hom(\ck,\ch)$.
(Note that $\psi_2\circ\gamma$ here is the $\psi_2$ in
the sum that we are estimating.)

What this shows is that the probability that we wish to estimate is
equal to the probability that a random $\theta\in\Theta$ is a
homomorphism from $\cd$ to $\ch$. Since we are assuming that $\ch$ is
$(\e,\cd,k-1)$-quasirandom and that $\e\le|\cd|^{-1}$, the counting
lemma (Corollary 5.2) implies that this is at most
$2\prod_{A\in\cd}\d_A=2\prod_{A\in\cd}\d_A$,
which proves the lemma. \hfill $\square$ \bigskip

Our next task is to show that we can make a small selection of
the functions $E_\psi$ and keep properties similar to those 
proved in the last two lemmas. The selection will be done in
the obvious way: randomly.
\medskip

\noindent {\bf Lemma {6.7}.} {\sl Let $\d=\prod_{A\in\cd}\d_A$,
let $\b=\prod_{A\in\ck}\d_A$ and let $r\ge\d^{-1}$ be a positive integer. 
Then there exist functions $\seq E r$ from $\Phi$ to $\{-1,0,1\}$ with the 
following three properties.}

(i) {\sl Each function $E_i$ is a product of $\{-1,0,1\}$-valued 
$A$-functions over subsets $A\subset B_0$ of size $k-1$.}

(ii) {\sl For each $i$ and each $\phi\in\Phi$, $E_i(\phi)$ is non-zero
only if $\phi(B_0)\in H_*([k])$.}

(iii) $\E_{i=1}^r\E_{\phi\in\Phi}F(\phi)E_i(\phi)\ge(\eta/2)\b$.

(iv) $\E_{\phi\in\Phi}\Bigl(\E_{i=1}^rE_i(\phi)\Bigr)^2\le 
(8\d/\eta\b)\E_{i=1}^r\E_{\phi\in\Phi}F(\phi)E_i(\phi)$.
\medskip

\Proof  For each $i$ let $E_i$ be one of the functions $E_\psi$, where
$\psi$ is chosen uniformly at random from $\Psi$. Let the choices be
independent (so, in particular, the $E_i$ are not necessarily distinct,
though they probably will be). Then it follows from Lemma 6.5 that
property (i) holds, and also that the expectation of 
$\E_{i=1}^r\E_{\phi\in\Phi}F(\phi)E_i(\phi)$ is at least $\eta\b$.

We now want to estimate the expectation of 
$\E_{\phi\in\Phi}\Bigl(\E_{i=1}^rE_i(\phi)\Bigr)^2$, and for this we
shall use Lemma 6.2, the technical lemma from the beginning of the section.
Set $n=|\Psi|=|\Phi|$ and let the vectors $\seq v n$ be the functions
$E_\psi$, which we regard as elements of $\L_2(\Phi)$. Lemma 6.6
tells us that $\bigl\|\E_{i=1}^rv_i\bigr\|_2^2\le 2\d$. Therefore,
Lemma 6.2 tells us that the expectation of 
$\bigl\|\E_{i=1}^rE_i\bigr\|_2^2$, which is the same as the expectation
of $\E_{\phi\in\Phi}\Bigl(\E_{i=1}^rE_i(\phi)\Bigr)^2$, is at most
$4\d$.

It follows that the expectation of
$$8\d\E_{i=1}^r\E_{\phi\in\Phi}F(\phi)E_i(\phi)-
\eta\beta\E_{\phi\in\Phi}\Bigl(\E_{i=1}^rE_i(\phi)\Bigr)^2$$
is at least $8\eta\b\d-4\eta\b\d=4\eta\b\d$.
It follows that there must be some choice of the functions $\seq E r$ 
such that the inequalities (iii) and (iv) are satisfied. 

Since each $E_i$ is one of the functions $E_\psi$, Lemma 6.5
implies that $E_i(\phi)$ is non-zero only if 
$(\phi,\psi)\in\Hom(\ck,\ch)$ for some $\psi\in\Psi$. But
a necessary condition for this is that $\phi(B_0)\in H_*([k])$,
so property (ii) is true as well. \hfill $\square$
\bigskip

\noindent {\bf Proof of Lemma 6.3.}  For each $i$ let us write
$E_i$ as a product $\prod_{j=1}^kE_{ij}$, where $E_{ij}$ is a
$\{-1,0,1\}$-valued $A_j$-function. (As in the proof of Lemma 6.5,
$A_j$ is the set $B_0\setminus\{(j,0)\}$.)

For each $j\le k$ we can partition the part $H_j$ of $H$ into at 
most $3^r$ sets, such that on each of these sets the function 
$E_{ij}$ is constant for every $i\le r$. Let $\seq Z N$ be the
corresponding induced partition of $H_*([k])$. (This concept 
was defined just before the statement of Lemma 6.3.) Then every
function $E_i$ is constant on every cell $Z_j$, from which it
follows that the function $g(\phi)=\E_{i=1}^rE_i(\phi)$ is constant
on every cell $Z_j$. (Here we are implicitly thinking of $g$
as a function of $\phi(B_0)$ and therefore defined on $H_*([k])$.)

With the help of Lemma 6.7, we are now in a position to apply Lemma
6.1. Property (iii) of Lemma 6.7 tells us that $\sp{F,g}\ge(\eta/2)\b$, 
and property (iv) tells us that
$\sp{F,g}/\|g\|_2^2\ge\eta\b/8\d$. 

Let $U$ be the set of all $\phi\in\Phi$ such that $\phi(B_0)\in
H_*([k])$. Then the map $\phi\mapsto\phi(B_0)$ is a bijection between $U$
and $H_*([k])$, so we can regard $\seq Z N$ as a partition of $U$, and
we can also regard $F$ and $g$ as functions defined on $U$. If we do
so, then their $L_2$-norms and inner products change: now we have
$\sp{F,g}\ge(\eta/2)\b/\zeta$, where $\zeta$ is the density of $U$ 
in $\Phi$, while the ratio $\sp{F,g}/\|g\|_2^2$ remains the same at 
$\ge\eta\b/8\d$.

Lemma 6.1 and these estimates tell us that the mean-square 
density of $F$ with respect to this partition of $U$ is at least
$(\eta\b/2\zeta)(\eta\b/8\d)=\eta^2\b^2/16\d\zeta$.
By Lemma 5.2 (the counting lemma), $\zeta\le 2\prod_{A\subsetneq B_0}\d_A$.
Recall that every set $A\subsetneq B_0$ is the index of precisely
$2^{|A|+1}-1$ sets in $\cd$ and $2^{|A|}$ sets in $\ck$. It
follows that $\b^2=\d\prod_{A\subsetneq B_0}\d_A\ge\d\zeta/2$. Therefore,
the mean-square density of $F$ with respect to the partition
$\seq Z N$ is at least $\eta^2/32$. Since $F(\phi)=f(\phi(B_0))$,
this statement is equivalent to the statement of 
Lemma~6.3.~\hfill~$\square$ \bigskip

\proclaim Corollary {6.8}. Let $\ch$ be a $k$-partite $(k-1)$-chain 
with vertex sets $\seq X k$, let $\cd$ be the double octahedron, let
$\d=\prod_{A\in\cd}\d_A$ and let $r\ge\d^{-1}$ be a positive integer.
Suppose that $\e\le|\cd|^{-1}$ and that $\ch$ is 
$(\e,\cd,k-1)$-quasirandom. Let $H^k$ be a $k$-partite $k$-uniform
hypergraph with vertex sets $\seq X k$, let the density of $H^k$ relative
to $\ch$ (that is, the quantity $|H^k|/|H_*([k])$) be $\d_{[k]}$ and 
suppose that $H^k$ is not $\eta$-quasirandom relative to $\ch$. Let $H$ be
the set of all edges of $\ch$ of size $k-1$ and let $\seq H k$ be the
$k$ parts of $H$. Then there are partitions of the $H_i$ into at most
$3^r$ sets each such that the mean-square density of (the characteristic
function of) $H^k$ with respect to the induced partition of $H_*([k])$ 
is at least $\d_{[k]}^2+\eta^2/32$.

\Proof  Let $f:H_*([k])\ra[-1,1]$ be the function $H^k-\d_{[k]}$. 
Then the statement that $H^k$ is not $\eta$-quasirandom relative to
$\ch$ is, by definition, the statement that $f$ is not
$\eta$-quasirandom relative to $\ch$. Therefore, by Lemma 6.3, we can
find partitions of the required kind for which the mean-square density
of $f$ with respect to the induced partition of $H_*([k])$ is at
least $\eta^2/32$. 

Let $\seq Z N$ be the induced partition of $H_*([k])$ and for each
$(\seq x k)\in Z_i$ let $G(\seq x k)=|H^k\cap Z_i|/|Z_i|$. Then
the mean of $G$ is the same as the mean of $H^k$, namely $\d_{[k]}$.
The value that $G$ takes in $Z_i$ can also be written as 
$\d_{[k]}+\E_{x\in Z_i}f(x)$, so the expectation of $(G-\d_{[k]})^2$,
which is also the mean-square density of $G-\d_{[k]}$ (since $G$ is
constant on the cells $Z_i$) is the mean-square density of $f$. But 
it is also the variance of $G$, so by the usual formula 
$\var X=\E X^2-(\E X)^2$ we find that the mean-square density of
$G$ is $\d_{[k]}^2$ plus the mean-square density of $f$. (Here 
we have again used the fact that $G$ is constant on cells, so that
the mean-square density of $G$ is just $\E G^2$.) The result 
follows.~\hfill~$\square$ \bigskip

\noindent {\bf \S 7. The statement of a regularity lemma for 
$r$-partite chains.}
\medskip

Corollary 6.8 is stage one of the proof of our regularity lemma. In 
this short section we will introduce some definitions and state
the regularity lemma itself. The proof (or rather, stage two of 
the proof) will be given in \S 9. 

Broadly speaking, the result says that we can take a $k$-uniform
hypergraph $H$, regard it as a chain (by adding all subsets of
edges of $H$) and decompose that chain into subchains almost all
of which are quasirandom. This is a useful thing to do, because
Corollary 5.2 gives us a good understanding of quasirandom chains.
Thus, the regularity lemma and counting lemma combine to allow
us to decompose any (dense) $k$-uniform hypergraph into pieces
that we can control. In the final section of the paper we shall
exploit this by proving a generalization of Theorems 1.3 and 1.6
to $k$-uniform hypergraphs, which implies the multidimensional
Szemer\'edi theorem.

Our principal aim will be to understand a certain $(k+1)$-partite
$k$-uniform hypergraph. However, for the purposes of formulating a
suitable inductive hypothesis it is helpful to prove a result that is
more general in two ways. First of all, we shall look at $r$-partite
$k$-uniform hypergraphs. Secondly, rather than looking at single
hypergraphs we shall look at {\it partitions}.  To be precise, let
$\seq X r$ be a sequence of finite sets. Given any subset
$A\subset[r]$, $A=\{\seq i s\}$, let $K(A)$ be the {\it complete
$s$-uniform hypergraph} on the sets $X_{i_1},\dots,X_{i_s}$, that is,
the hypergraph consisting of all subsets of $X_1\cup\dots\cup X_r$
that intersect $X_i$ in a singleton if $i\in A$ and are disjoint from
$X_i$ otherwise. For each $s\le r$, the {\it complete $s$-uniform
hypergraph} $K_s(\seq X r)$ on the sets $\seq X r$ is the union of the
hypergraphs $K(A)$ over all sets $A\subset[r]$ of size $s$.  Finally,
the {\it complete $k$-chain} on $\seq X r$, denoted $\ck_k(\seq X r)$,
is the union of all $K(A)$ such that $A$ has cardinality at most $k$:
that is, it consists of all subsets of $X_1\cup\dots\cup X_r$ of size
at most $k$ that intersect each $X_i$ at most once.

To form an arbitrary $r$-partite $s$-uniform hypergraph $H$ with 
vertex sets $\seq X r$, one can choose, for each $A\subset[r]$
of size $s$, a subset $H(A)\subset K(A)$ and let $H$ be the
union of these hypergraphs $H(A)$. If we want to, we can
regard each $H(A)$ as a partition of $K(A)$ into the two
sets $H(A)$ and $K(A)\setminus H(A)$. Our regularity lemma
will be concerned with more general partitions, but it will
imply a result for hypergraphs as an easy corollary.

Suppose now that for every subset $A\subset[r]$ of size 
at most $k$ we have a partition of the hypergraph $K(A)$.
If $B$ and $B'$ are two edges of this hypergraph (that is,
if they are two sets of index $A$), let us
write $B\sim_AB'$ if $B$ and $B'$ lie in the same cell
of the partition, and say that $B$ and $B'$ are $A$-{\it
equivalent}.

One can use these equivalence relations to define finer
ones as follows. Given two sets $B$, $B'$ of index
$A$ and given any subset $C\subset A$, there are unique
subsets $D\subset B$ and $D'\subset B'$ of index $C$.
Let us say that $B$ and $B'$ are $C$-{\it equivalent}
if $D$ and $D'$ are. Then let us say that $B$ and $B'$
are {\it strongly equivalent} if they are
$C$-equivalent for every subset $C\subset A$. In other words,
we ask not only for $B$ to belong to the same cell $B'$, but also
for every subset of $B$ to belong to the same cell as the
corresponding subset of $B'$ in the corresponding partition.

Given this system of equivalence relations, we can define
a collection of chains as follows. For every $r$-tuple
$x=(\seq x r)\in X_1\times\dots\times X_r$ and every set
$A$ of size at most $k$, let $x(A)$ be the set $\{x_i:i\in A\}$ 
and let $H(A,x)$ be the hypergraph consisting of all sets $B$ 
that are strongly equivalent to $x(A)$. 

\proclaim Lemma {7.1}. The union $\ch=\ch(x)$ of the hypergraphs 
$H(A,x)$ over all sets $A$ of size at most $k$ is an $r$-partite
$k$-chain.

\Proof  Let $B\in H(A,x)$ and let $D\subset B$. Let $C$ be the
index of $D$. Since $B$ is strongly equivalent to $x(A)$,
$D$ is strongly equivalent to $x(C)$. Therefore $D\in H(C,x)$
and the lemma is proved. \hfill $\square$ \bigskip

\proclaim Lemma {7.2}. Let $x=(\seq x r)$ and $y=(\seq y r)$
belong to the set $X_1\times\dots\times X_r$ and let $\ch(x)$
and $\ch(y)$ be the two chains constructed as in Lemma 7.1. Then 
for every set $A\subset[r]$ of size at most $k$, the hypergraphs
$H(A,x)$ and $H(A,y)$ are either equal or disjoint.

\Proof  Suppose that $B$ is a set of index $A$ and that
$B\in H(A,x)\cap H(A,y)$. Then $B$ is strongly equivalent
to both $x(A)$ and $y(A)$, so these two sets are strongly
equivalent to each other. It follows that $H(A,x)=H(A,y)$.
\hfill $\square$ \bigskip

Let us call two $r$-partite $k$-chains $\ch$ and $\ch'$ with the
same vertex sets $\seq X r$ {\it compatible} if, for every
subset $A\subset[r]$ of size at most $k$, the hypergraphs
$H(A)$ and $H'(A)$ are either equal or disjoint. By a
{\it chain decomposition} of the complete $r$-partite $k$-chain
$\ck_k(\seq X r)$ we mean a set $\{\seq \ch N\}$ of $r$-partite 
$k$-chains with the following two properties:

\itemitem{(i)} for every $i$ and $j$ the chains $\ch_i$ and $\ch_j$ are
compatible;

\itemitem{(ii)} for every sequence $x=(\seq x r)\in X_1\times\dots\times X_r$
there is precisely one chain from the set $\{\seq \ch N\}$ that
contains every subset of $\{\seq x r\}$ of size at most~$k$.
\smallskip

\noindent  Note that a chain decomposition is not a partition of 
$\ck_k(\seq X r)$. There is no interesting way to partition 
$\ck_k(\seq X r)$ into subchains, as a moment's thought will reveal. 
Lemmas 7.1 and 7.2 show that the chains $\ch(x)$ form a chain 
decomposition of $\ck_k(\seq X r)$. (It may be that $\ch(x)=\ch(y)$,
but this does not contradict (ii) because we have carefully defined 
a chain decomposition to be a {\it set} of chains rather than a 
{\it sequence} of chains.)

We are now ready to state our regularity lemma.

\proclaim Theorem {7.3}. Let $\cj$ be an
$r$-partite $k$-chain with vertex sets $\seq E r$ and let
$0<\e\le|\cj|^{-1}$. Let $\seq X r$ 
be a sequence of finite sets and for each subset $A\subset[r]$ 
of size at most $k$ let $\cp(A)$ be a partition of the 
hypergraph $K(A)$ into $n_A$ sets. Then there are refinements 
$\cq(A)$ of the partitions $\cp(A)$ leading to a chain decomposition
of $\ck_k(\seq X r)$ with the following property: if $x=(\seq x r)$
is a randomly chosen element of $X_1\times\dots\times X_r$ then 
the probability that the chain $\ch(x)$ is $(\e,\cj,k)$-quasirandom
is at least $1-\e$. Moreover, $\cq(A)=\cp(A)$ when $|A|=k$, and
for general $A$ the number of sets $m_A$ in the partition 
$\cq(A)$ depends only on $\e$, $\cj$, $k$ and the numbers 
$n_C$.

Before we start on the proof, let us comment on how we shall
actually use Theorem 7.3. We will be presented with an
$r$-partite $k$-uniform hypergraph $H$ with vertex sets
$\seq X r$. All the ${r\choose k}$ $k$-partite parts $H(A)$
of $H$ will have density at least a certain fixed $\d>0$.
We will then apply Theorem 7.3 to the partitions $\cp(A)$
defined as follows. If $|A|=k$ then $\cp(A)$ will be
$\{H(A),K(A)\setminus H(A)\}$. If $|A|<k$ then it will
be the trivial partition $\{K(A)\}$. In this case, the
result will tell us that we can find partitions $\cq(A)$
such that almost all edges of $\ch$ lie in quasirandom
chains from the decomposition determined by the partitions
$\cq(A)$.
\bigskip

\noindent {\bf \S 8. Basic facts about partitions and 
mean-square density.}
\medskip

In order to prove a regularity lemma for systems of partitions,
we need to generalize the notion of mean-square density as
follows. Let $\cp=\{\seq X r\}$ and $\cq=\{\seq Y s\}$ be two 
partitions of a finite set $U$. Then the 
{\it mean-square density of} $\cp$ {\it with respect to}
$\cq$ is the quantity
$$\sm i r\sm j s{|Y_j|\over |U|}\left({|X_i\cap Y_j|\over|Y_j|}\right)^2\ ,$$
that is, the sum of all the mean-square densities of the
sets $X_i$ (by which we mean the mean-square densities of their 
characteristic functions, as defined in \S 6) with respect 
to $\cq$.

Since the numbers $|X_i\cap Y_j|/|Y_j|$ are non-negative and
sum to 1, we have the simple upper bound
$$\sm i r\sm j s{|Y_j|\over |U|}\left({|X_i\cap Y_j|\over|Y_j|}\right)^2
\le\sm j s{|Y_j|\over |U|}=1$$
for this quantity. An alternative way of seeing this, which will be
helpful later, is to notice that each $u\in U$ is contained in a 
unique $X_i$ and a unique $Y_j$, and the mean-square density of
$\cp$ with respect to $\cq$ is the expected value of 
$|X_i\cap Y_j|/|Y_j|$.

\proclaim Lemma {8.1}. Let $\cp=\{\seq X r\}$ and $\cq=\{\seq Y s\}$ 
be two partitions of a finite set $U$, and let $\cq'$ be a refinement
of $\cq$. Then the mean-square density of $\cp$ with respect to
$\cq'$ is at least as great as the mean-square density of $\cp$
with respect to $\cq$.

\Proof  Let the sets that make up $\cq'$ be called $Y_{jk}$,
where $Y_j=\bigcup_kY_{jk}$. For each $j$ and $k$ define $\g_j$ 
and $\g_{jk}$ by $|Y_j|=\g_j|U|$ and $|Y_{jk}|=\g_{jk}|U|$.
For each $i$, $j$ and $k$ let $d_{ij}=|X_i\cap Y_j|/|Y_j|$
and let $d_{ijk}=|X_i\cap Y_{jk}|/|Y_{jk}|$. Then 
$$\sum_kd_{ijk}|Y_{jk}|=\sum_k|X_i\cap Y_{jk}|=|X_i\cap Y_j|=d_{ij}|Y_j|\ ,$$
from which it follows that $\sum_k\g_{jk}d_{ijk}=\g_jd_{ij}$ for
every $i$ and $j$.

The mean-square density of $\cp$ with respect to $\cq$ is
$\sum_i\sum_j\g_jd_{ij}^2$, which is therefore equal to
$$\eqalign{\sum_i\sum_j\g_j^{-1}\Bigl(\sum_k\g_{jk}d_{ijk}\Bigr)^2
&=\sum_i\sum_j\Bigl(\sum_k\g_j^{-1/2}\g_{jk}d_{ijk}\Bigr)^2\cr
&\le\sum_i\sum_j\Bigl(\sum_k\g_j^{-1}\g_{jk}\Bigr)
\Bigl(\sum_k\g_{jk}d_{ijk}^2\Bigr)\ ,\cr}$$
by the Cauchy-Schwarz inequality. Since $\sum_k\g_j^{-1}\g_{jk}=1$
for every $j$, this equals $\sum_i\sum_j\sum_k\g_{jk}d_{ijk}^2$,
which is the mean-square density of $\cp$ with respect to $\cq'$.
\hfill $\square$ \bigskip

The next lemma is a simple, but somewhat irritating, technicality.

\proclaim Lemma {8.2}. Let $\e>0$, let $\seq X r$ be a sequence 
of finite sets, let $\ck(\seq X r)$ be the complete $r$-partite
$k$-chain with vertex sets $\seq X r$ and for each
$A\subset\{1,2,\dots,r\}$ of size at most $k$ let $\cp(A)$ be a
partition of $K(A)$ into $n_A$ sets.  For each $x=(\seq x r)\in
X_1\times\dots\times X_r$ and each $A$ of size at most $k$ let
$\d_{A,x}$ be the relative density of the hypergraph $H(A,x)$ in the
chain $\ch(x)$ (defined in the previous section). Then if $(\seq x r)$
is chosen randomly from $X_1\times\dots\times X_r$ and
$A\subset\{1,2,\dots,r\}$ has size at most $k$, the probability that
$\d_{A,x}<\e n_A^{-1}$ is at most $\e$.

\Proof  Let $B$ and $B'$ be two sets of index $A$. Let us call
them {\it weakly equivalent}, and write $B\sim_*B'$, if
$B$ is $C$-equivalent to $B'$ for every {\it proper} subset $C$ 
of $A$. Then $B$ is strongly equivalent to $B'$ if and only if 
$B\sim_*B'$ and $B\sim_AB'$.

The relative density $\d_{A,x}$ is simply the probability that
a set $B$ of index $A$ is strongly equivalent to $x(A)$ given that
it is weakly equivalent to $x(A)$. Since $K(A)$ is partitioned into
$n_A$ sets, the number of strong equivalence classes in each
weak equivalence class is at most $n_A$. Therefore, for any
weak equivalence class $T$, the probability that $x(A)$ lies
in a strong equivalence class of size less than $\e n_A^{-1}|T|$ 
given that it lies in $T$ is at most $\e$. If $x(A)$ lies in a
strong equivalence class of size at least $\e n_A^{-1}|T|$, 
then the probability that $B$ is in the same strong equivalence
class given that $B$ is in $T$ is at least $\e n_A^{-1}$,
which implies that $\d_{A,x}\ge\e n_A^{-1}$.

Therefore, for every $T$ the conditional probability that 
$\d_{A,x}<\e n_A^{-1}$ given that $x(A)\in T$ is less than
$\e$. The result follows. \hfill $\square$ \bigskip

\noindent We now have all the ingredients needed to prove
our regularity lemma.
\bigskip

\noindent {\bf \S 9. The proof of Theorem 7.3.}  
\medskip

It will be convenient for the proof if for each set $A\subset[r]$ of
size at most $k$, the chain $\cj$ contains a copy $\cd_A$ of the
double octahedron of dimension $|A|$.  Since the result for $\cj$
follows from the result for any larger chain, we are free to assume that
this is the case.

We shall first describe an inductive procedure for producing better and 
better systems of partitions when the conclusion of Theorem 7.3 does not 
hold. Then we shall prove that the procedure terminates.

We shall need one piece of notation. Let $\seq X r$ be a sequence of
finite sets and for each subset $C\subset[r]$ of size at most $k$ let
$\cp(C)$ be a partition of the hypergraph $K(C)$. For each set $A\subset[r]$
of size at most $k$ we shall write $\sigma_A(\cp)$ for the mean-square
density of the partition $\cp(A)$ with respect to the partition of
$K(A)$ into weak equivalence classes with respect to the partition 
system $\cp$. (These were defined in the proof of Lemma 8.2 above.)
\bigskip

\noindent {\bf Lemma 9.1.} {\sl Let $\cj$ be an $r$-partite $k$-chain with 
vertex sets $\seq E r$ and let $0<\e\le |\cj|^{-1}$. Let $\seq X r$ be
a sequence of finite sets and for each subset $C\subset[r]$ of size at
most $k$ let $\cp(C)$ be a partition of the hypergraph $K(C)$ into
$n_C$ sets. For each $x=(\seq x r)$, let $\ch(x)$ be the chain arising
from $x$ and the corresponding chain decomposition of $\ck_k(\seq X
r)$. Suppose that when $x$ is chosen randomly from
$X_1\times\dots\times X_r$ the probability that $\ch(x)$ fails to be
$(\e,\cj,k)$-quasirandom is at least $\e$. Then there is a set $A$ of
size $s\le k$ and a system of refinements $\cq(C)$ of the partitions
$\cp(C)$ with the following properties.}

(i) {\sl $\cq(C)=\cp(C)$ and $\sigma_C(\cq)\ge\sigma_C(\cp)$ except 
if $C\subset A$ and $|C|=s-1$.}

(ii) {\sl $\sigma_A(\cq)$ exceeds $\sigma_A(\cp)$ by a non-zero amount
that depends only on $\cj$, $\e$, $k$ and the numbers of cells in 
the partitions $\cp(B)$ with $|B|\ge s$.}

(iii) {\sl When $C\subset A$ and $|C|=s-1$, the number of cells 
in the partition $C$ depends only on $\e$, $k$ and the numbers of 
cells in the partitions $\cp(B)$ with $B\subset C$.}
\bigskip

\Proof  For each set $C$, let $t_C$ be the number of
cells in the partition $\cp(C)$ of $K(C)$. 
Let $\g$ be defined by the equation $2\g\sm i k{r\choose i}=\e$.
By Lemma 8.2, the probability that there exists a subset
$C\subset[r]$ of size at most $k$ such that 
$\d_{C,x}<\g t_A^{-1}$ is at most $\g\sm i k{r\choose i}=\e/2$.
Therefore, with probability at least $\e/2$, the chain
$\ch(x)$ fails to be $(\e,\cj,k)$-quasirandom but for each
$C$ the relative density $\d_{C,x}$ is at least $\g t_C^{-1}$.

Let $\eta_2,\dots,\eta_k$ and $\e_2,\dots,\e_k$ be the sequences
that appear in the definition of quasirandom chains (in subsection
3.7), and note that $\eta_s$ depends only on
$\e$ and the densities $\d_{B,x}$ with $|B|\ge s$. Since 
$\d_{C,x}\ge\g t_C^{-1}$ for every $C$, it follows that $\eta_s$
is bounded below by a function of $\e$ and all those $t_B$ for
which $|B|\ge s$.

If $\ch(x)$ fails to be $(\e,\cj,k)$-quasirandom, then there must be 
a minimal $s$ such that it fails to be $(\e_s,\cj,s)$-quasirandom,
and for that $s$ there must be a set $A$ of size $s$ such that
$H(A,x)$ is not $\eta_s$-quasirandom relative to $\ch(x)$, while
$\ch(x)$ is $(\e_{s-1},\cj,s-1)$-quasirandom. Since there are
at most $\sm i k{r\choose i}$ possibilities for this set $A$
we may deduce from the last paragraph but one that there exists a
set $A$ of size $s\le k$ such that, with probability at
least $\g$, the chain $\ch(x)$ is $(\e_{s-1},\cj,s-1)$-quasirandom
but $H(A,x)$ is not $\eta_s$-quasirandom relative to $\ch(x)$ and
$\d_{C,x}\ge\g t_C^{-1}$ for every $C$.

Let us call $x$ {\it irregular} if $\ch(x)$ has these two properties.
Given an irregular $x$, let $\ch_-(A,x)$ be the $s$-partite
$(s-1)$-chain made up of all the hypergraphs $H(C,x)$ with
$C\subsetneq A$.  We can now apply Corollary 6.8 to the chain
$\ch_-(A,x)$ and to the $s$-uniform hypergraph $H(A,x)$. (Thus, the
$k$ of Corollary 6.8 is equal to $s$ here.) Since
$\e_{s-1}\le\e\le|\cj|^{-1}$ and $\cd_A\subset\cj$, the conditions
hold for the corollary to be applicable, with $k$ replaced by $s$. The
hypergraphs $\seq H k$ in the statement of Corollary 6.8 are, in this
context, the hypergraphs $H(A',x)$, where $A'$ ranges over all subsets
of $A$ of size $s-1$.

For each $C\subsetneq A$ we know that $\d_{C,x}\ge\g t_C^{-1}$. 
Therefore, if $r'$ is a positive integer that is at least 
$\prod_{C\in\cd_A}\g^{-1}t_C$, then for each subset $A'\subset A$
of size $s-1$ we can find a partition of $H(A',x)$ into at most $3^{r'}$
subsets, in such a way that the mean-square density of $H(A,x)$
with respect to the induced partition of $H_*(A,x)$ is at least
$\d_{A,x}^2+\eta_s^2/32$. (Here, $H_*(A,x)$ denotes the hypergraph consisting
of all sets $Y$ of index $A$ such that every proper subset of $Y$ 
belongs to $\ch_-(A,x)$.)

Let $\ch(A,x)$ be the $s$-partite $s$-chain $H(A,x)\cup\ch_-(A,x)$.
The number of distinct possibilities for $\ch(A,x)$ as $x$ varies is
at most $\prod_{C\subset A}t_C$.  For each one such that $x$ is
irregular (if $\ch(A,x)=\ch(A,y)$ and $x$ is irregular then $y$ is
irregular) choose a partition of the hypergraphs $H(A',x)$ as
above. In general, it will often happen that $\ch(A,x)\ne\ch(A,y)$ but
$H(A',x)=H(A',y)$, so each hypergraph $H(A',x)$ may be partitioned
many times. However, the number of distinct chains $\ch(A,x)$ is at
most $T_A=\prod_{C\subset A}t_C$, so we can find a common refinement
of all the partitions of $H(A',x)$ into at most $3^{r'T_A}$ sets.

For each $A'\subset A$ of size $s-1$ let $\cq(A')$ be the union of all
these common refinements, over all the different sets $H(A',x)$. There
are at most $T_{A'}$ of these sets, each partitioned into at most
$3^{r'T_A}$ sets, so $\cq(A')$ is a partition of $K(A')$ into at most
$T_{A'}3^{r'T_A}$ sets, and it refines the partition $\cp(A')$. For all
other sets $A$, let $\cq(A)=\cp(A)$.

By Lemma 8.1, given any irregular $x$, the mean-square density of 
$H(A,x)$ with respect to the partition of $H_*(A,x)$ that is 
induced by the refined partitions of the hypergraphs $H(A',x)$ is 
still at least $\d_{A,x}^2+\eta_s^2/32$. As for a regular $x$,
Lemma 8.1 tells us that the mean-square density of $H(A,x)$ 
with respect to the refined partition of $\ch(x)_*(A)$ is
still at least $\d_{A,x}^2$.

Let $\sigma_A(\cp)$ be the mean-square density of the partition $\cp(A)$
with respect to the partition of $K(A)$ into weak equivalence
classes coming from the partitions $\cp(C)$. Let $\sigma_A(\cq)$ be
the mean-square density of $\cp(A)=\cq(A)$ with respect to the
partition of $K(A)$ arising from $\cq$ in the same way. By the 
remark preceding Lemma 8.1, $\sigma_A(\cp)$ is the expectation of
$\d_{A,y}$ over all sequences $y=(\seq y r)$. Let us write this 
as $\d_{A,y}(\cp)$ since it depends on the system of partitions 
$\cp(C)$. Thus, $\sigma_A(\cp)$ is the expectation of 
$\d_{A,x}(\cp)$ and similarly for $\cq$.

What we have just shown is that if $x$ is irregular, then 
$\E[\d_{A,y}(\cq)|y\in H(A,x)]$ is at least 
$\d_{A,y}(\cp)^2+\eta_s^2/32$, which equals
$\E[\d_{A,y}(\cp)|y\in H(A,x)]+\eta_s^2/32$. If $x$ is regular,
then this conditional expectation is at least $\d_{A,y}(\cp)^2$, or
$\E[\d_{A,y}(\cp)|y\in H(A,x)]$. Since the probability that 
$x$ is irregular is at least $\g$, this shows that 
$\E[\d_{A,y}(\cq)]\ge\E[\d_{A,y}(\cp)]+\g\eta_s^2/32$. In
other words, $\sigma_A(\cq)\ge\sigma_A(\cp)+\g\eta_s^2/32$.

To summarize: if the conclusion of Theorem 7.3
is not true for the partitions $\cp(C)$ then there is a set
$A$ of size $s\le k$ and a system of refinements $\cq(C)$ such 
that $\cq(C)=\cp(C)$ except when $C$ is a subset of $A$ of size 
$s-1$, and such that $\sigma_A(\cq)\ge\sigma_A(\cp)+\g\eta_s^2/32$.
For a general $C$, we have $\sigma_C(\cq)\ge\sigma_C(\cp)$
except if $C\subset A$ and $|C|=s-1$. This is because if
$C$ is any other set, then $\cq(C)=\cp(C)$ and all other
partitions have either been refined or stayed the same.
Thus, the lemma is proved. \hfill $\square$ \bigskip

To complete the proof of Theorem 7.3, we must argue that this 
process of successive refinement cannot be iterated for ever.

Imagine, then, that we are trying to find an infinite sequence of
refinements of the kind we are given by Lemma 9.1. The difficulty we
face is that the mean-square densities $\sigma_C(\cp)$ tend to
increase, and there is always one set $A$ for which $\sigma_A(\cp)$
increases fairly substantially. Our only hope is that for the subsets
$C$ of $A$ obtained by removing one element, the mean-square densities
can drop considerably.

The trouble with that, however, is that the only way of getting the
mean-square density $\sigma_C(\cp)$ to drop is by getting the
mean-square density of some larger set $\sigma_A(\cp)$ to
increase. 

To see why this observation leads to a proof, suppose that we do
indeed have an infinite sequence of refinements of the kind given to
us by Lemma 9.1. Then there must be a set $A$ of maximal cardinality
$s$ that is used infinitely many times. It follows that there must be
some point in the sequence after which $A$ is used infinitely many
times but no set of larger cardinality is ever used. After that point,
the only partitions $\cp(C)$ that change are for sets $B$ of
cardinality less than $s$, by (i) of Lemma 9.1. It follows from (ii)
that after that point the quantity $\sigma_A(\cq)$ increases
infinitely often by an amount that does not change as the iteration
proceeds. This is a contradiction, since $\sigma_A(\cq)$ is bounded
above by~1.  The proof of the regularity lemma is complete.
\bigskip

A careful examination of the above argument shows that the bound that
arises from it increases by one level in the Ackermann hierarchy each
time $k$ increases by 1, except at the jump from the trivial case
$k=1$ to the first non-trivial case $k=2$, when we go from nothing to a
bound of tower type. In particular, since we shall need $k$-uniform
hypergraphs to prove the multidimensional Szemer\'edi theorem for sets
of size $k+1$, our bound for that theorem is of Ackermann type. The
only cases where better bounds are known are the one-dimensional case,
which is treated in [G1], and the case of sets of size 3, where a
trebly exponential bound was obtained by Shkredov [S].
\bigskip

\noindent {\bf \S 10. Hypergraphs with few simplices.}
\medskip

Now that we have established counting and regularity lemmas we have
the tools necessary to prove the generalization of Theorems 1.3 and
1.6 to $k$-uniform hypergraphs.  

\proclaim Theorem {10.1}. Let $k$ be a positive integer. Then for
every $a>0$ there exists $c>0$ with the following property. Let 
$H$ be a $(k+1)$-partite $k$-uniform hypergraph with vertex sets 
$\seq X {{k+1}}$, and let $N_i$ be the size of $X_i$. Suppose that
$H$ contains at most $c\prod_{i=1}^{k+1}N_i$ simplices. Then for 
each $i\le k+1$ one can remove at most $a\prod_{j\ne i}N_j$ 
edges of $H$ from $\prod_{j\ne i}X_j$ in such a way that after 
the removals one is left with a hypergraph that is simplex-free.

\Proof  For each subset $A\subset[k+1]$ of size at most $k$, define 
a partition $\cp(A)$ of $K(A)$ as follows. If $|A|<k$ then $\cp(A)$ 
consists of the single set $K(A)$. If $|A|=k$ then it consists of
the sets $H(A)$ and $K(A)\setminus H(A)$. Now apply Theorem 7.3
to this system of partitions, with $\cj=[k+1]^{(\le k)}$ and
$\e=\min\{|\cj|^{-1}/2,a/2\}$, obtaining for each $A\in\cj$ a partition
$\cq(A)$ of $K(A)$ into $m_A$ sets.

If $x=(\seq x {{k+1}})\in X_1\times\dots\times X_{k+1}$ and $\ch(x)$
is not $(\e,\cj,k)$-quasirandom, then there must be some $A$ of 
size $s\le k$ such that $H(A,x)$ is not $\eta_s$-quasirandom relative
to $\ch(x)$. There must be some $i$ such that $i\notin A$, and if
$(\seq y {{k+1}})$ is another sequence such that $y_j=x_j$ when
$j\ne i$, then $H(A,y)$ will also not be $\eta_s$-quasirandom
relative to $\ch(y)$. Therefore, since $\ch(x)$ is 
$(\e,\cj,k)$-quasirandom with probability at least $1-\e$, there 
are at most $\e\prod_{j\ne i}N_j$ elements of $\prod_{j\ne i}X_j$
that can be extended to sequences $x$ such that $\ch(x)$ is not 
$(\e,\cj,k)$-quasirandom. Remove from $H$ any such element.

Let $\g$ be defined by $\g\sm i k{k+1\choose i}=a/2$. Lemma 8.2 
tells us that if $x=(\seq x {{k+1}})$ is chosen randomly, then
with probability at least $1-a/2$, we have $\d_{A,x}\ge\g m_A^{-1}$
for every $A\in[k+1]^{(\le k)}$. Again, the event that this happens
for a particular $A$ does not depend on the $x_i$ with $i\notin A$.
So for each $i$ there are at most $a\prod_{j\ne i}N_j/2$ elements 
of $\prod_{j\ne i}X_j$ that can be extended to sequences $x$ for
which $\d_{A,x}<\g m_A^{-1}$ for some $A\subset[k+1]$ with $i\notin A$.
Once again, remove all such elements from $H$.

For each $i$ we have removed at most $a\prod_{j\ne i}N_j$ elements
from $H\cap\prod_{j\ne i}X_j$. It remains to show that in the process
we have either removed all simplices from $H$, or else, for some
$c>0$ that depends on $a$ only, there were at least $c\prod_jN_j$ 
simplices to start with.

Suppose, then, that after the removals there is still a simplex 
$x=(\seq x {{k+1}})$, and consider the chain $\ch(x)$. Then
for every $A\subset[k+1]$ of size $k$ the following statements
are true. First, the set $x(A)$ is an element of $H$ (or else
$x$ would not be a simplex). Second, the hypergraph $H(A,x)$
is a subset of $H$ (since $x(A)\in H$ and the partition 
into strong equivalence classes resulting from $\cq$ refines 
the partition $\cp$). Third, $\d_{C,x}\ge\g m_C^{-1}$
for every $C\subset A$ (or else we would have removed $x(A)$
from $H$). Finally, the chain $\ch(x)$ is $(\e,\cj,k)$-quasirandom
(or else for some $A$ of size $k$ we would have removed $x(A)$
from $H$).

We now apply Corollary 5.2, the counting lemma for quasirandom
chains. It implies that the number of simplices in the chain
$\ch(x)$, which is the same as the number of homomorphisms from
$\cj$ to $\ch(x)$, is at least $\prod_jN_j\prod_{A\in\cj}\d_{A,x}$,
which is at least $\prod_jN_j\prod_{A\in\cj}\g m_A^{-1}$. But
$\g$ and the $m_A$ depend on $a$ and $k$ only, so the result
is proved. \hfill $\square$ \bigskip

Finally, let us deduce from this a multidimensional Szemer\'edi
theorem.

\proclaim Theorem {10.2}. Let $\d>0$ and $k\in\N$. Then, if $N$
is sufficiently large, every subset $A$ of the $k$-dimensional
grid $\{1,2,\dots,N\}^k$ of size at least $\d N^k$ contains a
set of points of the form $\{a\}\cup\{a+de_i:1\le i\le k\}$,
where $\seq e k$ is the standard basis of $\R^k$ and $d$ is
a non-zero integer.

\Proof  Suppose that $A$ is a subset of $\{1,2,\dots,N\}^k$
of size $\d N^k$, and that $A$ contains no configuration of the kind
claimed. Define a $(k+1)$-partite $k$-graph $F_k$ with vertex sets
$\seq X {{k+1}}$ as follows. If $j\le k$ then the elements of $X_j$
are hyperplanes of the form $P_{j,m}=\{(\seq x k):x_j=m\}$ for some
integer $m\in\{1,2,\dots,N\}$.  If $j=k+1$ then they are hyperplanes
of the form $Q_m=\{(\seq x k):x_1+\dots+x_k=m\}$ where $m$ is an
integer between $k$ and $kN$. The edges of $F_k$ are sets of $k$
hyperplanes from different sets $X_j$ that intersect in a point of
$A$.

If $F_k$ contains a simplex with vertices $P_{j,m_j}$ and $Q_m$,
then the points $(\seq m k)$ and $(\seq m k)+(m-\sm i km_i)e_j$
all belong to $A$. This gives us a configuration of the desired
kind except in the degenerate case where $m=\sm i km_i$, which
is the case where all $k+1$ hyperplanes have a common intersection. 
By our assumption on $A$, all the simplices in $F_k$ are therefore
degenerate ones of this kind, which implies that there are at
most $\d N^k$ of them.

Now $|X_i|=N$ if $i\le k$ and $|X_{k+1}|=kN$. We can therefore apply
(the contrapositive of) Theorem 10.1 with $c=N^{-1}k^{-1}$.  If $N$ is
sufficiently large, then the resulting $a$ is smaller than $\d/2k$,
which implies that we can remove fewer than $\d N^k$ edges from the
hypergraph $F_k$ and thereby remove all simplices. However, every edge
of a degenerate simplex determines the point of intersection of the
$k+1$ hyperplanes and hence the simplex itself. It follows that one
must remove at least $\d N^k$ edges to get rid of all simplices. This
contradiction proves the theorem. \hfill $\square$ \bigskip

The above result is a special case of the multidimensional
Szemer\'edi theorem, but it is in fact equivalent to the
whole theorem. This is a well-known observation. We give
a (slightly sketchy) proof below.

\proclaim Theorem {10.3}. For every $\d>0$, every positive
integer $r$ and every finite subset $X\subset\Z^r$ there
is a positive integer $N$ such that every subset $A$ of
the grid $\{1,2,\dots,N\}^r$ of size at least $\d N^r$
has a subset of the form $a+dX$ for some positive integer
$d$.

\Proof  It is clearly enough to prove the result for sets
$X$ such that $X=-X$, so all we actually need to ensure is that $d\ne
0$. A simple averaging argument shows that we may also assume that $X$
is not contained in any $(r-1)$-dimensional subspace of $\R^r$. Let
the cardinality of $X$ be $k+1$. Let $\phi$ be an affine map that
defines a bijection from the set $\{0,\seq e k\}\subset\R^k$ to $X$,
regarded as a subset of $\R^r$. Another simple averaging argument
allows us to find a grid $\{1,2,\dots,M\}^k$, where $M$ tends to
infinity with $N$, as well as a point $z\in\Z^r$ and a constant 
$\eta>0$ depending on $\d$ and $X$ only, such that $z+\phi(x)\in A$
for at least $\eta M^k$ points in $\{1,2,\dots,M\}^k$. Let $B$ be
the set of points with this property. Thus, $B$ has density at least
$\eta$ and Theorem 10.2 shows that $B$ contains a set of the form 
$w+c\{0,\seq e k\}$. But then $z+\phi(w+c\{0,\seq e k\})$ is a
set of the form $a+dX$ and is also a subset of $A$. \hfill $\square$
\bigskip

\noindent {\bf Concluding Remarks.} 
\medskip

This paper has a slightly strange history, which may be worth briefly
outlining here. The main results were first obtained in 2003, and a
preprint circulated. I am very grateful indeed to Yoshiyasu Ishigami,
who read this preprint carefully and found an error which, though it
did not invalidate the approach, occurred early in the argument and
therefore necessitated changes throughout the paper. While thinking
about how to go about this rewriting, I discovered a much simpler
proof of the counting lemma, and in the end it seemed best, even if
depressing, to rewrite the whole paper (including the regularity part)
from scratch.

I owe a second debt of gratitude to the two referees, who also read
the paper with great care. Not only did they save me from a large
number of minor errors, but they also made valuable suggestions about
the presentation of the paper. While thinking about how to respond to
these suggestions I realized, with a certain sense of d\'ej\`a vu,
that the sections on the counting lemma could still be greatly
improved. The argument that now appears is essentially the same, but
the notation has been changed and the triple induction slightly
reorganized, with the result that the proof is now shorter, clearer, and
easier to identify with the arguments presented in the special cases
in \S 2. That section, as was mentioned in the footnote at the
beginning of it, was not in the original version of the paper. The
excellent idea of presenting some small examples was suggested by 
one of the referees.

In 2005, Tao [T] gave another proof of the main result of this paper
(Theorem 10.1), and indeed of a slight generalization. He too proved
regularity and counting lemmas. His methods were more closely related
to those of Nagle, R\"odl, Schacht and Skokan, but he introduced some
new ideas and a different language that led to considerably shorter
proofs than theirs.
\bigskip

\noindent {\bf References.}
\medskip

\noindent [FK] H. Furstenberg, Y. Katznelson, {\it An ergodic
Szemer\'edi theorem for commuting transformations}, J. Analyse Math.
{\bf 34} (1978), 275-291.
\medskip

\noindent [FKO] H. Furstenberg, Y. Katznelson, D. Ornstein, 
{\it The ergodic theoretical proof of Szemer\'edi's theorem}, 
Bull. Amer. Math. Soc. {\bf 7} (1982), 527--552.
\medskip

\noindent [G1] W. T. Gowers, {\it A new proof of Szemer\'edi's theorem},
Geom. Funct. Anal. {\bf 11} (2001), 465-588.
\medskip

\noindent [G2] W. T. Gowers, {\it Quasirandomness, Counting and 
Regularity for 3-Uniform Hypergraphs}, Combin. Probab. Comput. {\bf 15}
(2006), 143-184.
\medskip

\noindent [NRS] B. Nagle, V. R\"odl and M. Schacht, {\it The counting
lemma for regular $k$-uniform hypergraphs}, Random Structures and Algorithms
{\bf 28} (2006), 113-179.
\medskip

\noindent [R] V. R\"odl, {\it Some developments in Ramsey theory}, 
Proceedings of the International Congress of Mathematicians, Vol. I, II
(Kyoto 1990), 1455-1466, Math. Soc. Japan, Tokyo, 1991.
\medskip

\noindent [RS] V. R\"odl and J. Skokan, {\it Regularity lemma for 
$k$-uniform hypergraphs}, Random Structures and Algorithms {\bf 25} (2004), 
1-42.

\noindent [Ro] K. Roth, {\it On certain sets of integers}, J. London Math.
Soc. {\bf 28} (1953), 245-252.
\medskip

\noindent [RS] I. Z. Ruzsa and E. Szemer\'edi, {\it Triple systems with
no six points carrying three triangles}, Combinatorics (Proc. Fifth
Hungarian Colloq., Keszthely, 1976), Vol. II, 939-945.
\medskip

\noindent [S] Shkredov, I. D., {\it On a problem of Gowers}, (Russian)  
Izv. Ross. Akad. Nauk Ser. Mat. {\bf 70}  (2006), 179--221;  translation 
in  Izv. Math. {\bf 70} (2006),  385--425
\medskip

\noindent [So1] J. Solymosi, {\it Note on a generalization of Roth's theorem},
Discrete and Computational Geometry, 825-827, Algorithms Combin. {\bf 25},
Springer, Berlin 2003.
\medskip

\noindent [So2] J. Solymosi, {\it A note on a question of Erd\H os and
Graham}, Combin. Probab. Comput. {\bf 13} (2004), 263-267.
\medskip

\noindent [Sz1] E. Szemer\'edi, {\it Integer sets containing no 
$k$ elements in arithmetic progression}, Acta Arith. {\bf 27} (1975),
299-345.
\medskip

\noindent [Sz2] E. Szemer\'edi, {\it Regular partitions of graphs}, in
Probl\`emes Combinatoires et Th\'eorie des Graphes, Proc. Colloque
Inter. CNRS, (Bermond, Fournier, Las Vergnas, Sotteau, eds.),
CNRS Paris, 1978, 399-401.
\medskip

\noindent [T] T. Tao, {\it A variant of the hypergraph removal lemma},
J. Combin. Theory Ser. A {\bf 113} (2006), 1257-1280.

\bigskip

\bye